\documentclass[11pt]{amsart}
\usepackage[utf8]{inputenc}
\usepackage{enumerate,amssymb,xcolor,comment}
\usepackage[a4paper, footskip=1cm, headheight = 16pt, top=3.4cm, bottom=2.5cm,  right=2.5cm,  left=2.5cm]{geometry}
\usepackage[colorlinks,linkcolor=blue,citecolor=red]{hyperref}
\usepackage[center]{caption}

\newtheorem{theorem}{Theorem}[section]
\newtheorem{lemma}[theorem]{Lemma}

\newtheorem{corollary}[theorem]{Corollary}

\usepackage[foot]{amsaddr}

\usepackage[T1]{fontenc}

\usepackage[leqno]{amsmath}

\makeatletter
\newcommand{\leqnomode}{\tagsleft@true}
\newcommand{\reqnomode}{\tagsleft@false}
\makeatother

\usepackage{latexsym}
\usepackage{amsfonts}

\usepackage{tikz}
\usetikzlibrary{shapes}
\usepackage{float}
\usepackage{lmodern}

\def\dd{\hbox{-}}

\usepackage{marvosym}               

\DeclareMathOperator{\tw}{tw}
\DeclareMathOperator{\width}{width}

\newcounter{tbox}

\newcommand{\sta}[1]{\medskip\refstepcounter{tbox}\noindent{ \parbox{\textwidth}{(\thetbox) \emph{#1}}}\vspace*{0.3cm}}

\newcommand{\mylongtitle}[1]{%
  \ifodd\value{page}%
    \protect\parbox{0.97\linewidth}{#1}\hfill%
  \else%
    \hfill\protect\parbox{0.97\linewidth}{#1}%
  \fi%
}

\usepackage{latexsym}
\usepackage{amsfonts}
\usepackage[leqno]{amsmath}
\usepackage{tikz}
\usepackage{float}
\usepackage{lmodern}
\usepackage[T1]{fontenc}

\def\dd{\hbox{-}}

\usepackage{marvosym}

\newcommand{\otherlabel}[2]{\protected@edef\@currentlabel{#2}\label{#1}}
\makeatother

\mathchardef\mh="2D

\title[Induced subgraphs and tree decompositions VII.]{Induced subgraphs and tree decompositions\\
VII. Basic obstructions in $H$-free graphs\footnote{This is an accepted manuscript. The published version appeared in the Journal of Combinatorial Theory, Series B
Volume 164, January 2024, Pages 443--472 and is available here: \url{https://doi.org/10.1016/j.jctb.2023.10.008}}}

\author{Tara Abrishami$^{\ast \dagger}$}
\author{Bogdan Alecu$^{\ast \ast \mathparagraph}$}
\author{Maria Chudnovsky$^{\ast \amalg}$}
\author{Sepehr Hajebi $^{\mathsection}$}
\author{Sophie Spirkl$^{\mathsection \parallel}$}

\address{$^{\ast}$Princeton University, Princeton, NJ, USA}
\address{$^{**}$School of Computing, University of Leeds, Leeds, UK}
\address{$^{\mathsection}$Department of Combinatorics and Optimization, University of Waterloo, Waterloo, Ontario, Canada}
\address{$^{\dagger}$ Supported by NSF-EPSRC Grant DMS-2120644.}
\address{$^{\amalg}$ Supported by NSF-EPSRC Grant DMS-2120644 and by AFOSR grant FA9550-22-1-0083.} 
     \address{$^{\mathparagraph}$ Supported by DMS-EPSRC Grant EP/V002813/1.} 
\address{$^{\parallel}$ We acknowledge the support of the Natural Sciences and Engineering Research Council of Canada (NSERC), [funding reference number RGPIN-2020-03912].
Cette recherche a \'et\'e financ\'ee par le Conseil de recherches en sciences naturelles et en g\'enie du Canada (CRSNG), [num\'ero de r\'ef\'erence RGPIN-2020-03912]. This project was funded in part by the Government of Ontario.}

\date {\today}

\allowdisplaybreaks

\begin{document}
\maketitle

\begin{abstract} 
We say a class $\mathcal{C}$ of graphs is {\em clean} if for every positive integer $t$ there exists a positive integer $w(t)$ such that every graph in $\mathcal{C}$ with treewidth more than $w(t)$ contains an induced subgraph isomorphic to one of the following: the complete graph $K_t$, the complete bipartite graph $K_{t,t}$, a subdivision of the $(t\times t)$-wall or the
line graph of a subdivision of the $(t \times t)$-wall. In this paper, we adapt a method due to Lozin and Razgon (building on earlier ideas of Wei\ss auer) to prove that the class of all \textit{$H$-free} graphs (that is, graphs with no induced subgraph isomorphic to a fixed graph $H$) is clean if and only if $H$ is a forest whose components are subdivided stars.

Their method is readily applied to yield the above characterization. However, our main result is much stronger: for every forest $H$ as above, we show that forbidding certain connected graphs containing $H$ as an induced subgraph (rather than $H$ itself) is enough to obtain a clean class of graphs. Along the proof of the latter strengthening, we build on a result of Davies and produce, for every positive integer $\eta$, a complete description of unavoidable connected induced subgraphs of a connected graph $G$ containing $\eta$ vertices from a suitably large given set of vertices in $G$. This is of independent interest, and will be used in subsequent papers in this series.

\end{abstract}

\section{Introduction}\label{intro}

\subsection*{Background} All graphs in this paper are finite and simple.

Treewidth is a well-studied graph parameter that is of great interest in both structural and algorithmic graph theory. It was notably featured in the seminal work of Robertson and Seymour on graph minors \cite{RS-GMII}, and in numerous other papers ever since. For a more in-depth overview of the literature, the reader is invited to see, for example, Bodlaender's survey \cite{Bodlaendersurvey} and the references therein. 

\begin{figure}[t]
    \centering
    \includegraphics[scale=0.9]{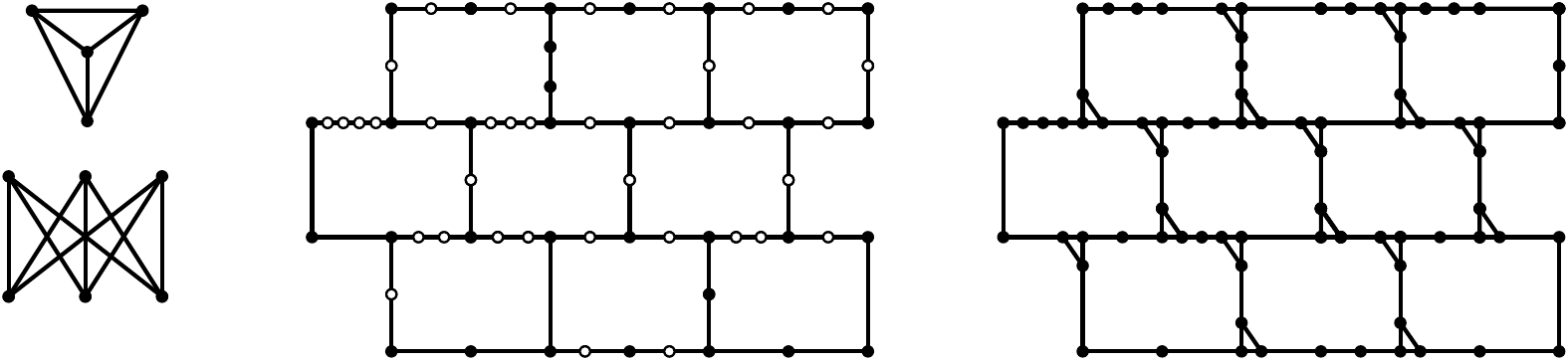}
    \caption{The $3$-basic obstructions, including a subdivision of $W_{3\times 3}$ (middle) and its line graph (right).}
    \label{fig:3basic}
\end{figure}

As a part of their graph minors series, Robertson and Seymour fully described the unavoidable minors in graphs of large treewidth. The relevant result, the so-called Grid Theorem \cite{RS-GMV}, states that every graph of large enough treewidth must contain a minor isomorphic to a large grid, or equivalently, a subgraph isomorphic to a large wall (the {\em $(t\times t)$-wall}, denoted by $W_{t\times t}$, is a planar graph of maximum degree three on $2t^2-2t$ vertices; see \cite{wallpaper} for a precise definition and see Figure~\ref{fig:3basic}). Since walls have large treewidth themselves, and treewidth cannot increase when taking minors, that result gives a structural dichotomy: a graph has large treewidth if and only if it contains a large wall as a subgraph.

The overarching goal of the current series and of several other recent works \cite{aboulker, multistar, Korhonen, lozin, layered-wheels, longhole2} is to understand treewidth from the perspective of induced subgraphs rather than minors. A first remark is that to force bounded treewidth, we need to forbid four kinds of induced subgraphs: a complete graph $K_t$, a complete bipartite graph $K_{t, t}$, all subdivisions of the $(t\times t)$-wall $W_{t\times t}$ for some $t$, and the line graphs of all subdivisions of $W_{t\times t}$ for some $t$. Let us call these graphs the {\em $t$-basic obstructions} (see Figure~\ref{fig:3basic}), and say that a graph $G$ is {\em $t$-clean} if $G$ contains no induced subgraph isomorphic to a $t$-basic obstruction. Moreover, we say a class $\mathcal{C}$ of graphs is \textit{clean} if the treewidth of $t$-clean graphs in $\mathcal C$ is bounded from above by a function of $t$.

The class of all graphs is not clean: various constructions of unbounded treewidth avoiding the basic obstructions have been discovered \cite{multistar, daviesconst, layered-wheels}. In fact, it is at the moment unclear whether a dichotomy similar to the Grid Theorem is at all achievable for induced subgraphs. Nevertheless, steady progress is being made. Of note is the following result, characterizing all finite sets of graphs which yield bounded treewidth when forbidden as induced subgraphs: 

\begin{theorem}[Lozin and Razgon \cite{lozin}]\label{lozinmain}
Let $\mathcal{H}$ be a finite set of graphs. Then the class of all graphs with no induced subgraph isomorphic to a member of $\mathcal{H}$ has bounded treewidth if and only if $\mathcal{H}$ contains a complete graph, a complete bipartite graph, a forest of maximum degree at most three in which every component has at most one vertex of degree more than two, and the line graph of such a forest.
\end{theorem}

In addition, several clean classes have been identified. For instance, Aboulker, Adler, Kim, Sintiari
and Trotignon \cite{aboulker} proved that every proper minor-closed class of graphs is clean:

\begin{theorem}[Aboulker, Adler, Kim, Sintiari
and Trotignon \cite{aboulker}]\label{minorclosed}
For every graph $H$, the class of all graphs with no minor isomorphic to $H$ is clean. Equivalently, for every graph $H$ and integers $t\geq 1$, there exists an integer $\xi=\xi(H,t)\geq 1$ such that every graph with no minor isomorphic to $H$ and treewidth more than $\xi$ contains either a subdivision of $W_{t\times t}$ or the line graph of a subdivision of $W_{t\times t}$ as an induced subgraph.
\end{theorem}

They also conjectured that graph classes of bounded maximum degree are clean, which was later proved by Korhonen \cite{Korhonen}:

\begin{theorem}[Korhonen \cite{Korhonen}]\label{boundeddeg}
For every integer $d\geq 1$, the class of graphs of maximum degree at most $d$ is clean. Equivalently, for all integers $d,t\geq 1$, there exists an integer $\gamma=\gamma(d,t)\geq 1$ such that every graph with maximum degree at most $d$ and treewidth more than $\gamma$ contains either a subdivision of $W_{t\times t}$ or the line graph of a subdivision of $W_{t\times t}$ as an induced subgraph.
\end{theorem}

There are also a number of results concerning holes, where a \textit{hole} in a graph is an induced cycle of length at least four. In particular, it was shown that (even hole, diamond, pyramid)-free graphs are clean \cite{pyramiddiamond}, and graphs in which no vertex has two or more neighbors in a hole disjoint from itself are clean \cite{onenbr}. It was also independently proved twice that graphs with no long hole are clean. For every positive integer $\lambda$, let $\mathcal{H}_{\lambda}$ be the class of all graphs with no hole of length more than $\lambda$.

\begin{theorem}[Gartland, Lokshtanov, Pilipczuk, Pilipczuk and Rz\k{a}\.{z}ewski \cite{longholes}, Wei\ss auer \cite{longhole2}]\label{longholepili}
For every integer $\lambda\geq 1$, the class $\mathcal{H}_{\lambda}$ is clean.
\end{theorem}

\subsection*{Our results}

The main result of this paper is Theorem~\ref{mainconnectify}. The precise statement of Theorem ~\ref{mainconnectify} requires some set-up, and we postpone it until Section~\ref{mainresultssec}. Informally, we show that every $t$-clean graph of sufficiently large treewidth contains, as an induced subgraph, a ``connectification'' of any given subdivided star forest $F$. Roughly speaking, this is a graph which can be partitioned into a ``rooted'' copy of $F$ and a second part, which only attaches at the roots of $F$ and ``minimally connects'' these roots. 

The proof of Theorem~\ref{mainconnectify} uses three ingredients. The first one is Theorem~\ref{blockmanystar}, which adapts the methods from \cite{lozin} (itself employing the strategy from \cite{longhole2}) in order to show that clean graphs with a large {\em block} -- a certain kind of highly connected structure -- must contain a large subdivided star forest. As a byproduct of this, we also obtain another way to derive Theorem~\ref{longholepili}. 

The second ingredient is Theorem~\ref{noblock}. This theorem combines a result of Wei\ss auer linking blocks and tree decompositions, together with Korhonen's bounded degree result (Theorem~\ref{boundeddeg}), in order to show that the class of graphs without a large block is clean.

The final ingredient, Theorem~\ref{minimalconnectedgeneral}, is a result of independent interest, and will be used in future papers in our series. Starting from a result of Davies \cite{Davies}, we provide a complete description of minimal connected graphs containing many vertices from a suitably large subset of a connected component. Put differently, we show that if a large enough set of vertices belongs to the same component, then a large subset of them are contained in one of a few prescribed induced subgraphs.

\medskip

We note that the first two out of those intermediate results already yield (the difficult direction of) an appealing dichotomy for clean classes defined by one forbidden induced subgraph. Indeed, writing $\mathcal{F}_H$ for the class of graphs with no induced subgraph isomorphic to $H$, we prove:

\begin{theorem}\label{mainstarforest1}
Let $H$ be a graph. Then $\mathcal{F}_H$ is clean if and only if $H$ is a subdivided star forest.
\end{theorem}

While the stronger Theorem~\ref{mainconnectify} might appear unwieldy at first, we remark that it has easier-to-state implications that are still more general than the above dichotomy. To illustrate this, denote by $\tilde{\mathcal{F}}_H$ the class of all graphs with no induced subgraph isomorphic to a subdivision of $H$. It follows that the ``if'' direction of Theorem~\ref{mainstarforest1} is equivalent to $\tilde{\mathcal{F}}_H$ being clean for every subdivided star forest $H$, and Theorem~\ref{longholepili} is equivalent to $\tilde{\mathcal{F}}_H$ being clean for every cycle $H$. Then Theorem~\ref{mainconnectify} readily implies the following, where by a \textit{subdivided double star}, we mean a a tree with at most two vertices of degree more than two.

\begin{theorem}\label{mainsubforest}
Let $H$ be a forest in which one component is a subdivided double star and every other component is a subdivided star. Then $\tilde{\mathcal{F}}_H$ is clean.
\end{theorem}
\begin{figure}
    \centering
    \includegraphics[scale=0.4]{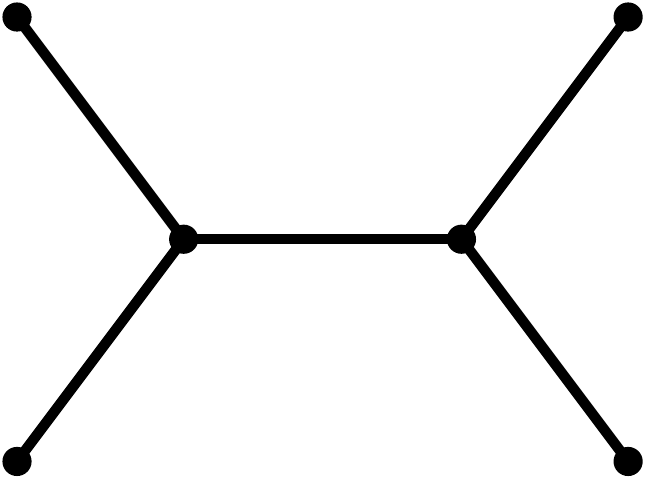}
    \caption{The smallest tree that is not a subdivided star.}
    \label{fig:doublestar}
\end{figure}
We remark that a full grid-type theorem for induced subgraphs is equivalent to a characterization of families $\mathcal{H}$ of graphs for which the class of all $\mathcal{H}$-free graphs is clean. This remains out of reach, and Theorem~\ref{mainstarforest1} takes the first step towards answering this question by characterizing all singletons $\mathcal{H}$ for which the class of all $\mathcal{H}$-free graphs is clean.

Here is a natural next step: \textit{for which \textbf{finite} families $\mathcal{H}$ of graphs is the class of all $\mathcal{H}$-free graphs clean?} From Theorem~\ref{mainstarforest1}, it follows that such a finite set $\mathcal{H}$ containing a subdivided star forest has the above property. One may then speculate that in fact \textit{all} finite set of graphs with the above property must contain a subdivided star forest. This, however, is false: for instance, assume that $H$ is the unique double star on six vertices (see Figure~\ref{fig:doublestar}; note that $H$ is the smallest tree that is not a subdivided star). Then $\mathcal{H}=\{H,K_3\}$ has the above property; in fact, this follows from the main result of an upcoming paper \cite{tw13} where the last four authors of the present work provide a full description of finite families $\mathcal{H}$ for which the class of all $\mathcal{H}$-free graphs is clean.

\subsection*{Outline of the paper} 
 We set up our notation and terminology in Section~\ref{defns}. Section~\ref{daviessec} describes the construction of \cite{daviesconst}, which is used to prove the ``only if'' direction of Theorem~\ref{mainstarforest1}. In Section~\ref{mainresultssec}, we state Theorem~\ref{mainconnectify} precisely, and show how to deduce Theorems~\ref{mainstarforest1} and~\ref{mainsubforest} from it. In Section~\ref{connectedsec}, we show that a connected graph $G$ with a sufficiently large subset $S$ of its vertices contains an induced connectifier with many vertices from $S$. The main result of Section~\ref{noblocksec} is Theorem~\ref{noblock}, where we prove that the class of graphs with no $k$-block is clean. In Section~\ref{distancesec}, we show that in a $t$-clean graph, every huge block can be transformed into a large block such that there is no short path between any two vertices of the new block. Section~\ref{ramseyblocksec} uses this in order to show that a $t$-clean graph with a huge block contains a large subdivided star forest. Finally, in Section~\ref{connectifysec}, we combine the main results from Sections~\ref{connectedsec}, \ref{noblocksec} and \ref{ramseyblocksec} to prove Theorem~\ref{mainconnectify}.

\section{Preliminaries}
\label{defns}
\subsection*{Graphs, subgraphs, and induced subgraphs} All graphs in this paper are finite and with no loops or multiple edges. Let $G=(V(G),E(G))$ be a graph. A \textit{subgraph} of $G$ is a graph obtained from $G$ by removing vertices or edges, and an \textit{induced subgraph} of $G$ is a graph obtained from $G$ by only removing vertices. Given a subset $X \subseteq V(G)$, $G[X]$ denotes the subgraph of $G$ induced by $X$, that is, the graph obtained from $G$ by removing the vertices not in $X$. We put $G \setminus X = G[V(G) \setminus X]$ (and in general, we will abuse notation and use induced subgraphs and their vertex sets interchangeably). Additionally,  for an edge $e\in E(G)$, we write $G-e$ to denote the graph obtained from $G$ by removing the edge $e$. For a graph $H$, by a \textit{copy of $H$ in $G$}, we mean an induced subgraph of $G$ isomorphic to $H$, and we say \textit{$G$ contains $H$} if $G$ contains a copy of $H$. We also say $G$ is \emph{$H$-free} if $G$ does not contain $H$.  For a class $\mathcal{H}$ of graphs we say $G$ is $\mathcal{H}$-free if $G$ is $H$-free for every $H \in \mathcal{H}$. For a graph $H$, we write $G=H$ whenever $G$ and $H$ have the same vertex set and the same edge set. 

\subsection*{Neighborhoods} Let $v \in V(G)$. The \emph{open neighborhood of $v$}, denoted by $N(v)$, is the set of all vertices in $G$ adjacent to $v$. The \emph{closed neighborhood of $v$}, denoted by $N[v]$, is $N(v) \cup \{v\}$. Let $X \subseteq G$. The \emph{open neighborhood of $X$}, denoted by $N(X)$, is the set of all vertices in $G \setminus X$ with at least one neighbor in $X$.
If $H$ is an induced subgraph of $G$ and $X \subseteq G$ with $H\cap X=\emptyset$, then $N_H(X)=N(X) \cap H$.
Let $X,Y \subseteq V(G)$ be disjoint. We say $X$ is \textit{complete} to $Y$ if all possible edges with one end in $X$ and one end in $Y$ are present in $G$, and $X$ is \emph{anticomplete}
to $Y$ if there is no edge between $X$ and $Y$. In the case $X=\{x\}$, we often say $x$ is \textit{complete} (\textit{anticomplete}) to $Y$ to mean $X$ is complete (anticomplete) to $Y$.

\subsection*{Tree decompositions and blocks} A \emph{tree decomposition} $(T, \chi)$ of $G$ consists of a tree $T$ and a map $\chi: V(T) \to 2^{V(G)}$ with the following properties: 
\begin{enumerate}[(i)]
\itemsep -.2em
    \item For every vertex $v \in V(G)$, there exists $t \in V(T)$ such that $v \in \chi(t)$. 
    
    \item For every edge $v_1v_2 \in E(G)$, there exists $t \in V(T)$ such that $v_1, v_2 \in \chi(t)$.
    
    \item For every vertex $v \in V(G)$, the subgraph of $T$ induced by $\{t \in V(T) \mid v \in \chi(t)\}$ is connected.
\end{enumerate}

For each $t\in V(T)$, we refer to $\chi(t)$ as a \textit{bag of} $(T, \chi)$.  The \emph{width} of a tree decomposition $(T, \chi)$, denoted by $\width(T, \chi)$, is $\max_{t \in V(T)} |\chi(t)|-1$. The \emph{treewidth} of $G$, denoted by $\tw(G)$, is the minimum width of a tree decomposition of $G$. 

\subsection*{Cliques, stable sets, paths, and cycles}
A \textit{clique} in $G$ is a set of pairwise adjacent vertices in $G$, and a \textit{stable set} in $G$ is a set of pairwise non-adjacent vertices in $G$.  A {\em path in $G$} is an induced subgraph of $G$ that is a path, while a {\em cycle in $G$} is a (not necessarily induced) subgraph of $G$ that is a cycle.  If $P$ is a path, we write $P = p_1 \dd \cdots \dd p_k$ to mean that $V(P) = \{p_1, \dots, p_k\}$, and $p_i$ is adjacent to $p_j$ if and only if $|i-j| = 1$. We call the vertices $p_1$ and $p_k$ the \emph{ends of $P$}, and say that $P$ is \emph{from $p_1$ to $p_k$}. The \emph{interior of $P$}, denoted by $P^*$, is the set $P \setminus \{p_1, p_k\}$. For a path $P$ in $G$ and $x,y\in P$, we denote by $P[x,y]$ the subpath of $P$ with ends $x$ and $y$. The \emph{length} of a path $P$ is the number of its edges. Let $C$ be a cycle. We write $C = c_1 \dd \cdots \dd c_k\dd c_1$ to mean $V(C) = \{c_1, \dots, c_k\}$, and $c_i$ is adjacent to $c_j$ if and only if $|i-j|\in \{1,k-1\}$. A {\em hole of $G$} is an induced subgraph of $G$ that is a cycle. The {\em length} of a cycle or a hole is the number of its edges.

\subsection*{Subdivisions} By a \textit{subdivision} of a graph $G$, we mean a graph obtained from $G$ by replacing the edges of $G$ by pairwise internally disjoint paths between the corresponding ends. Let $r\geq 0$ be an integer. An $r$-\textit{subdivision} of $G$ is a subdivision of $G$ in which the path replacing each edge has length $r+1$. Also, a $(\leq r)$-\textit{subdivision} of $G$ is a subdivision of $G$ in which the path replacing each edge has length at most  $r+1$, and a $(\geq r)$-\textit{subdivision} of $G$ is defined similarly. We refer to a $(\geq 1)$-subdivision of $G$ as a \textit{proper} subdivision of $G$.

\subsection*{Classes of graphs} A class $\mathcal{C}$ of graphs is called \textit{hereditary} if it is closed under isomorphism and taking induced subgraphs, or equivalently, if $\mathcal{C}$ is the class of all $\mathcal{H}$-free graphs for some family $\mathcal{H}$ of graphs. For a class of graphs $\mathcal{C}$ and a positive integer $t$, we denote by $\mathcal{C}^t$ the class of all $t$-clean graphs in $\mathcal{C}$. Thus, $\mathcal{C}$ is clean if for every positive integer $t$ there exists a positive integer $w(t)$ such that every graph in $\mathcal{C}^t$ has treewidth at most $w(t)$. The following is immediate from the definition of a clean class. 
\begin{lemma}\label{cleanlemma}
Let $\mathcal{X}$ be a class of graphs. Assume that for every $t$, there exists a clean class of graphs $\mathcal{Y}_t$ such that $\mathcal{X}^t\subseteq \mathcal{Y}_t$. Then $\mathcal{X}$ is clean. In particular, every subclass of a clean class is clean.
\end{lemma}

\subsection*{Forests and stars} By a \textit{branch vertex} of a graph $G$, we mean a vertex of degree more than two in $G$. For every forest $F$, we say a vertex $v\in V(F)$ is a \textit{leaf} of $F$ if $v$ has degree at most one in $F$. We denote by $\mathcal{L}(F)$ the set of all leaves of $F$.  By a \textit{star} we mean a graph isomorphic to the complete bipartite graph $K_{1,\delta}$ for some integer $\delta\geq 0$, and a \textit{star forest} is a forest in which every component is a star. Then subdivided stars are exactly trees with at most one branch vertex, and subdivided star forests are exactly forests in which every component is a subdivided star. A \textit{subdivided double star} is a tree with at most two branch vertices. 

By a \textit{rooted subdivided star} $S$ we mean a subdivided star $S$ together with a choice of one vertex $r$ in $S$, called the \textit{root}, such that if $S$ is not a path, then $r$ is the unique branch vertex of $S$. \textit{A rooted subdivided star forest} $F$ is a subdivided star forest with a choice of a root for every component of $F$. We also refer to the root of each component of $F$ as a \textit{root} of $F$, and denote by $\mathcal{R}(F)$ the set of all roots of $F$. By a \textit{stem} in $F$, we mean a path in $F$ from a leaf to a root. It follows that each stem is the (unique) path from a leaf of some component of $F$ to the root of the same component. The \textit{reach} of a rooted subdivided star $S$ is the maximum length of a stem in $S$. Also, the \textit{reach} of a subdivided star forest $F$ is the maximum reach of its components and the \textit{size} of $F$ is the number of its components. For a positive integer $\theta$ and graph $H$, we denote by $\theta H$ the disjoint union of $\theta$ copies of $H$. For integers $\delta\geq 0$ and $\lambda\geq 1$, we denote by $S_{\delta,\lambda}$ the $(\lambda-1)$-subdivision of $K_{1,\delta}$. So for $\delta\geq 3$, $\theta S_{\delta,\lambda}$ is a subdivided star forest of maximum degree $\delta$, reach $\lambda$ and size $\theta$.

\section{A construction from \cite{daviesconst}} \label{daviessec}
The goal of this section is to prove the ``only if'' direction of Theorem~\ref{mainstarforest1} using a construction from \cite{daviesconst}.

We begin with a definition, which will be used in subsequent sections, as well. Let $P$ be a path and $\rho,\sigma\geq 0$ and $\theta\geq 1$ be integers. A $2\theta$-tuple $(p_1,\ldots,p_{2\theta})$ of vertices of $P$ is said to be a $(\rho,\sigma)$-\textit{widening} of $P$ if
\begin{itemize}
    \item the vertices $p_1$ and $p_{2\theta}$ are the ends of $P$;
    \item traversing $P$ from $p_1$ to $p_{2\theta}$, the vertices $p_1,\ldots,p_{2\theta}$ appear on $P$ in this order;
    \item $P[p_{2i-1}, p_{2i}]$ has length $\rho$ for each $i\in [\theta]$, and;
    \item $P[p_{2i}, p_{2i+1}]$ has length at least $\sigma$ for each $i\in [\theta-1]$.
\end{itemize}

The $(\rho,\sigma)$-widening $(p_1,\ldots,p_{2\theta})$ is \textit{strict} if for each $i\in [\theta-1]$, $P[p_{2i}, p_{2i+1}]$ has length equal to $\sigma$. Also, we say a $\theta$-tuple $(p_1,\ldots,p_{\theta})$ of vertices of $P$ is a $\sigma$-\textit{widening} of $P$ if the $2\theta$-tuple $(p_1,p_1\ldots,p_{\theta},p_{\theta})$ is a $(0,\sigma)$-\textit{widening} of $P$.

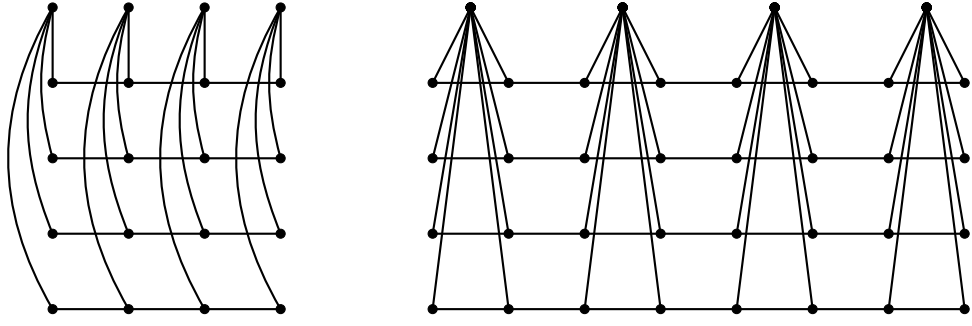
\begin{figure}[t]

\begin{tikzpicture}[scale=1,auto=left]
\tikzstyle{every node}=[inner sep=1.2pt, fill=black,circle,draw]  

        \foreach \i in {1,...,4} {
        \draw [thick] (-1, \i) -- (-4 ,\i);
        \node() at (-\i, 5) {};
        
        \draw [thick] (-\i , 5) edge[out = -90, in = 90] (-\i , 4);
        
        \draw [thick] (-\i , 5) edge[out = -105, in = 105] (-\i  , 3);
        
        \draw [thick] (-\i , 5) edge[out = -112, in = 112] (-\i , 2);
        
        \draw [thick] (-\i , 5) edge[out = -120, in = 120] (-\i , 1);

    }

\foreach \i in {1,...,4} {
            \foreach \x in {1,...,4} {
            \node() at (-\i, \x) {};
                
            }
        }

        \foreach \i in {1,...,4} {
        \draw [thick] (1, \i) -- (8 ,\i);
        
        \draw [thick] (\i+\i-0.5 , 5) edge[] (\i+\i , 4);
        \draw [thick] (\i+\i-0.5 , 5) edge[]  (\i+\i-1 , 4);
 
 \draw [thick] (\i+\i-0.5 , 5) edge[] (\i+\i , 3);
        \draw [thick] (\i+\i-0.5 , 5) edge[] (\i+\i-1 , 3);

         \draw [thick] (\i+\i-0.5 , 5) edge[] (\i+\i , 2);
        \draw [thick] (\i+\i-0.5 , 5) edge[] (\i+\i-1 , 2);

         \draw [thick] (\i+\i-0.5 , 5) edge (\i+\i , 1);
        \draw [thick] (\i+\i-0.5 , 5) edge (\i+\i-1 , 1);

    }

\foreach \i in {1,...,8} {
            \foreach \x in {1,...,4} {
            \node() at (\i,\x) {};
\node() at (\x+\x-0.5, 5) {};

            }
        }

\end{tikzpicture}
\caption{The graphs $J_{0,1,4}$ (left) and $J_{1,1,4}$ (right).}
\label{daviesfig}
\end{figure}

We now describe the construction of \cite{daviesconst} (though \cite{daviesconst} only mentions the case $\rho=0$). Let $\rho\geq 0$, $\sigma\geq 1$ and $\theta\geq 2$ be integers. We define $J=J_{\rho,\sigma,\theta}$ to be the graph with the following specifications (see Figure~\ref{daviesfig}).
\begin{itemize}
    \item $J$ contains $\theta$ pairwise disjoint and anticomplete paths $P_1,\ldots, P_{\theta}$.
    \item For each $j\in [\theta]$, $P_j$ admits a strict $(\rho,\sigma)$-widening $(p^j_1,\ldots,p^j_{2\theta})$.
    \item We have $J\setminus (\bigcup_{i\in [\theta]}V(P_i))=\{x_1,\ldots,x_{\theta}\}$ such that $x_1,\ldots,x_{\theta}$ are all distinct, and for all $i,j\in [\theta]$, we have $N_{J}(x_i)=\bigcup_{j\in [\theta]}P_j[p^j_{2i-1}, p^j_{2i}]$.
\end{itemize}
 
The following was proved in \cite{Davies}. Here we include a proof for the sake of completeness.

\begin{theorem}\label{daviesclean}
    For all integers $\rho\geq 0$, $\sigma\geq 1$ and $\theta\geq 2$, $J_{\rho,\sigma,\theta}$ is a $4$-clean graph of treewidth at least $\theta$.
\end{theorem}
\begin{proof}
     Note that $J_{\rho,\sigma,\theta}$ contains a $K_{\theta,\theta}$-minor (by contracting each path $P_i$ into a vertex), which implies that $\tw(J_{\rho,\sigma,\theta})\geq \theta$. Also, $J_{\rho,\sigma,\theta}$ is easily seen to be $\{K_4,K_{3,3}\}$-free. Let us say that a connected graph $H$ is \textit{feeble} if either $H$ has a vertex $v$ such that $H\setminus N_H[v]$ is not connected, or $H$ has a set $S$ of at most two branch vertices such that $H\setminus S$ has maximum degree at most two. Then every connected induced subgraph of $J_{\rho,\sigma,\theta}$ is feeble. On the other hand, for an integer $t\geq 4$, let $H$ be either a subdivision of $W_{t\times t}$ or the line graph of such a subdivision. Then one may observe that for every vertex $v\in H$, $H\setminus N_H[v]$ is connected. Moreover, $H$ contains a stable set $S$ of branch vertices with $|S|\geq 3$. It follows that $H$ is not feeble, and so $H$ is not isomorphic to an induced subgraph of $J_{\rho,\sigma,\theta}$. Hence, $J_{\rho,\sigma,\theta}$ is $4$-clean, as desired.
\end{proof}

The proof of the next lemma is straightforward, and we leave it to the reader.

\begin{lemma}\label{daviesgirth}
    For all integers $\sigma\geq 1$ and $\theta\geq 2$, the following hold.
    \begin{itemize}
        \item $J_{0,\sigma,\theta}$ has girth at least $2\sigma+4$.
        \item Let $u_1,u_2\in J_{1,\sigma,\theta}$ such that for each $i\in \{1,2\}$, $N_{J_{1,\sigma,\theta}}(u_i)$ contains a stable set of cardinality three. Then there is no path of length less than $\sigma+2$ in $J_{1,\sigma,\theta}$ from $u_1$ to $u_2$.
    \end{itemize}
\end{lemma}

We are now ready to prove the main result of this section.

\begin{theorem}\label{onlyifmain}
    Let $H$ be a graph for which $\mathcal{F}_H$ is clean. Then $H$ is a subdivided star forest.
\end{theorem}
\begin{proof}
    By the assumption, for every integer $t\geq 1$, there exists an integer $w(t)\geq 1$ such that every $t$-clean graph in $\mathcal{F}_H$ has treewidth at most $w(t)$. We deduce:

    \sta{\label{daviesforest}$H$ is a forest.}

    Suppose not. Let $\sigma$ be the length of the shortest cycle in $H$. By Theorem~\ref{daviesclean}, $J_{0,\sigma,w(4)+1}$ is $4$-clean. Also, by the first outcome of Lemma~\ref{daviesgirth}, $J_{0,\sigma,w(4)+1}$ has girth at least $2\sigma+4$, and so $J_{0,\sigma,w(4)+1}\in \mathcal{F}_H$. But then we have $\tw(J_{0,\sigma,w(4)+1})\leq w(4)$, which violates Theorem~\ref{daviesclean}. This proves \eqref{daviesforest}.

      \sta{\label{daviesstarforest}Every component of $H$ has at most one branch vertex.}

      Suppose for a contradiction that some component $C$ of $H$ contains two branch vertices $u$ and $v$. By \eqref{daviesforest}, $H$ is a forest, and so $C$ is a tree. Therefore, there exists a unique path in $H$ from $u$ to $v$, say of length $\sigma$, and we have $|N_{H}(u)\setminus N_{H}(v)|,|N_{H}(v)\setminus N_{H}(u)|\geq 2$. It follows from the second outcome of Lemma~\ref{daviesgirth} that $J_{1,\sigma,w(4)+1}\in \mathcal{F}_H$. Also, by Theorem~\ref{daviesclean}, $J_{1,\sigma,w(4)+1}$ is $4$-clean. But then we have $\tw(J_{1,\sigma,w(4)+1})\leq w(4)$, a contradiction with Theorem~\ref{daviesclean}. This proves \eqref{daviesstarforest}.\medskip

      Now the result follows from \eqref{daviesforest} and \eqref{daviesstarforest}. This completes the proof of Theorem~\ref{onlyifmain}.
\end{proof}

\section{Connectification and statement of the main result}
\label{mainresultssec}
Here we state the main result of the paper, Theorem~\ref{mainconnectify}. Then we discuss how it implies Theorems~\ref{mainstarforest1} and \ref{mainsubforest}.

We need numerous definitions. A vertex $v$ of a graph $G$ is said to be \textit{simplicial} if $N_G(v)$ is a clique of $G$. The set of all simplicial vertices of $G$ is denoted by $\mathcal{Z}(G)$. It follows that every degree-one vertex in $G$ belongs to $\mathcal{Z}(G)$. In particular, for every forest $F$, we have $\mathcal{L}(F)=\mathcal{Z}(F)$.

By a \textit{caterpillar} we mean a tree $C$ of maximum degree three in which all branch vertices lie on a path. A path $P$ in $C$ is called a \textit{spine} for $C$ if all branch vertices of $C$ belong to $V(P)$ and subject to this property $P$ is maximal with respect to inclusion (our definition of a caterpillar is non-standard for two reasons: a caterpillar is often allowed to be of arbitrary maximum degree, and a spine often contains all vertices of degree more than one.)

Let $C$ be a caterpillar with $\theta\geq 3$ leaves. Note that $C$ has exactly $\theta-2$ branch vertices, and both ends of each spine of $C$ are leaves of $C$. Also, for every leaf $l\in \mathcal{L}(C)$, there exists a unique branch vertex in $C$, denoted by $v_l$, for which the unique path in $C$ from $l$ to $v_l$ does not contain any branch vertex of $C$ other than $v_l$ (and, in fact, $\{v_l:l\in \mathcal{L}(C)\}$ is the set of all branch vertices of $C$). We say an enumeration $(l_1,\ldots, l_{\theta})$ of $\mathcal{L}(C)=\mathcal{Z}(C)$ is $\sigma$-\textit{wide} if for some spine $P$ of $C$, the $\theta$-tuple $(l_1, v_{l_2},\ldots ,v_{l_{\theta-1}}, l_{\theta})$ is a $\sigma$-widening of $P$. Also, let $H$ be the line graph of $C$. Then assuming $e_l$ to be the unique edge in $C$ incident with the leaf $l\in \mathcal{L}(C)$, we have $\mathcal{Z}(H)=\{e_l:l\in \mathcal{L}(C)\}$. An enumeration $(e_{l_1},\ldots, e_{l_{\theta}})$ of $\mathcal{Z}(H)$ is called $\sigma$-\textit{wide} if $(l_1,\ldots, l_{\theta})$ is a $\sigma$-wide enumeration of $\mathcal{L}(C)$.  By a $\sigma$-\textit{caterpillar}, we mean a caterpillar $C$ for which $\mathcal{L}(C)$ admits a $\sigma$-wide enumeration. It follows that if $H$ is the line graph of a caterpillar $C$, then $\mathcal{Z}(H)$ admits a $\sigma$-wide enumeration if and only if $C$ is a $\sigma$-caterpillar.

Let $H$ be a graph and $S$ be a set. We say $H$ is \textit{$S$-tied} if $\mathcal{Z}(H)\subseteq H\cap S$ and \textit{loosely $S$-tied} if $\mathcal{Z}(H)=H\cap S$. Also, for a positive integer $\eta\geq 1$, we say $H$ is (\textit{loosely}) $(S,\eta)$-\textit{tied} if $H$ is (loosely) $S$-tied and $|H\cap S|=
\eta$. It follows that if $H$ is loosely $(S,\eta)$-tied, then $|\mathcal{Z}(H)| = \eta$.

For a graph $G$, a set  $S\subseteq G$ and integers $\eta\geq 2$ and $\sigma\geq 1$ and $i\in \{0,\ldots, 4\}$, we say an induced subgraph $H$ of $G$ is an $(S,\eta,\sigma)$-\textit{connectifier of type $i$} if $H$ satisfies the condition (C$i$) below.

\begin{enumerate}[(C1)]
 \setcounter{enumi}{-1}
   \item $H$ is a loosely $(S,\eta)$-tied line graph of a subdivided star in which every stem has length at least $\sigma$.
  \item $H$ is an $(S,\eta)$-tied rooted subdivided star with root $r$ in which every stem has length at least $\sigma$, and we have $(H\cap S)\setminus \mathcal{L}(H)\subseteq \{r\}$. 
  \item $H$ is an $(S,\eta)$-tied path with $H\cap S=\{s_1,\ldots, s_{\eta}\}$ where $(s_1,\ldots, s_{\eta})$ is a $\sigma$-widening of $H$.
  \item $H$ is a loosely $(S,\eta)$-tied $\sigma$-caterpillar.
  \item $H$ is a loosely $(S,\eta)$-tied line graph of a $\sigma$-caterpillar.
  \end{enumerate}

\begin{figure}[t]
    \centering
    \includegraphics[scale=0.8]{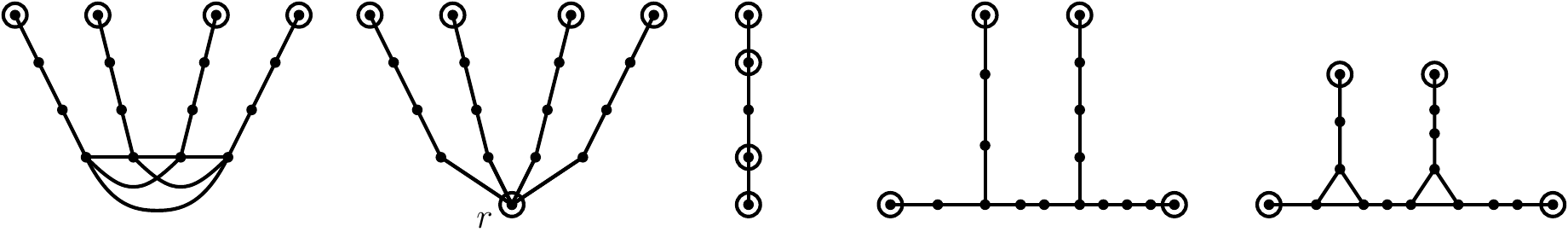}
    \caption{From left to right: an $(S,4)$-connectifier $H$ of type $0,1,2,3$ and $4$. Circled nodes depict the vertices in $H\cap S$. Note that for the subdivided star, $r$ may or may nor belong to $S$ (and if it does, then we have $\eta=5$).}
    \label{fig:connectifiers}
\end{figure}

See Figure~\ref{fig:connectifiers}. We say $H$ is an $(S,\eta)$-\textit{connectifier of type $i$} if it is an $(S,\eta,1)$-connectifier of type $i$. Also, we say $H$ is an $(S,\eta,\sigma)$-\textit{connectifier} (resp.\ $(S,\eta)$-\textit{connectifier}) if it is an $(S,\eta,\sigma)$-connectifier (resp.\ $(S,\eta)$-connectifier) of type $i$ for some $i\in \{0,\ldots, 4\}$. 

Note that connectifiers of type $0$ contain large cliques, and since we mostly work with $t$-clean graphs, they do not come up in our arguments. However, for the sake of generality, we cover them in both the above definition and the main result of 
 the next section, Theorem~\ref{minimalconnectedgeneral}. We also remark that, unlike the connectifiers of other types, connectifiers of type $1$ in fact need to be ``tied'' rather than ``loosely tied.'' For instance, let $G$ be a subdivided star with root $r$ and let $S=\mathcal{L}(G)\cup \{r\}$. Then for every $\eta>1$, every $(S,\eta)$-connectifier in $G$ contains $r$.

Let $\sigma$ be a positive integer,  $F$ be a graph and $X\subseteq F$ with $|X|\geq 2$. Let $\pi:[|X|]\rightarrow X$ be a bijection. By a $\sigma$-\textit{connectification of $(F,X)$ with respect to $\pi$}, we mean a graph $\Xi$ with the following specifications.
\begin{itemize}
\item $F$ is an induced subgraph of $\Xi$.
\item $F\setminus X$ is anticomplete to $\Xi\setminus F$.
\item Let $H=\Xi\setminus (V(F)\setminus X)$. Then $H$ is $(X,|X|,\sigma)$-connectifier in $\Xi$ of type $i$ for $i\in [4]$ such that
\begin{itemize}
 \item if $H$ is of type 2 (that is, $H$ is path), then, traversing $H$ from one end to another, $(\pi(1),\ldots, \pi(|X|))$ is a $\sigma$-widening of $H$, and;
     \item if $H$ is of type $3$ or $4$, then $(\pi(1),\ldots, \pi(|X|))$ is a $\sigma$-wide enumeration of $\mathcal{Z}(H)$.
\end{itemize}
\end{itemize}

Also, by a $\sigma$-\textit{connectification of $(F,X)$}, we mean a $\sigma$-connectification of $(F,X)$ with respect to some bijection $\pi:[|X|]\rightarrow X$.

Let $\mathcal{C}_{\sigma, F,X,\pi}$ be the class of all graphs with no induced subgraph isomorphic to a $\sigma$-connectification of $(F,X)$ with respect to $\pi$, and $\mathcal{C}_{\sigma,F,X}$ be the class of all graphs with no induced subgraph isomorphic to a $\sigma$-connectification of $(F,X)$. In other words, $\mathcal{C}_{\sigma, F,X}$ is the intersection of all classes $\mathcal{C}_{\sigma,F,X, \pi}$ over all bijections $\pi:[|X|]\rightarrow X$. As a result, for every $\pi:[|X|]\rightarrow X$, we have $\mathcal{C}_{\sigma, F,X}\subseteq \mathcal{C}_{\sigma,F,X, \pi}$.\\

The following is our main result, which we will prove in Section~\ref{connectifysec}.

\begin{theorem}\label{mainconnectify}
    Let $\sigma\geq 1$ be an integer, $F$ be a rooted subdivided star forest of size at least two and $\pi:[|\mathcal{R}(F)|]\rightarrow \mathcal{R}(F)$ be a bijection. Then the class $\mathcal{C}_{\sigma, F, \mathcal{R}(F),\pi}$ is clean. 
\end{theorem}

Next we discuss briefly how to deduce Theorems~\ref{mainstarforest1} and \ref{mainsubforest} using Theorem~\ref{mainconnectify}. The ``only if'' direction of Theorem~\ref{mainstarforest1} is proved in Theorem~\ref{onlyifmain}. Also, the ``if'' direction of Theorem~\ref{mainstarforest1} follows from Theorem~\ref{mainsubforest}. So it suffices to prove Theorem~\ref{mainsubforest}, which we restate:

\begin{theorem}\label{mainsubforestrestate}
Let $H$ be a forest in which one component is a subdivided double star and every other component is a subdivided star. Then $\tilde{\mathcal{F}}_H$ is clean.
\end{theorem}

\begin{proof}
We define $F$ and  $\sigma$ as follows. If $H$ is a subdivided star forest, then let $F=2H$ be rooted and $\sigma=2$. If $H$ is not a subdivided star forest, let $H'$ be the $1$-subdivision of $H$. Then there are two branch vertices $u_1,u_2\in H'$ and a path $Q$ in $H'$ from $u_1$ to $u_2$ with $Q^*\neq \emptyset$ such that $F'=H'\setminus Q^*$ is a subdivided star forest. For each $i\in \{1,2\}$, let $F_i$ be the component of $F'$ containing $u_i$. Then $u_i$ is a vertex of maximum degree in $F_i$ and so $u_i$ is a valid choice for a root of $F_i$. Let $F'$ be rooted such that $u_1,u_2\in \mathcal{R}(F')$. Let $\delta, \lambda$ and $\theta$ be the maximum degree, the reach and the size of $F'$, respectively. So we have $\delta,\theta\geq 2$ and $\lambda\geq 1$. Let $F=\theta S_{\delta+1,\lambda}$ be rooted with its unique choice of roots and let $\sigma=|Q|\geq 3$. Then, every $\sigma$-connectification of $(F,\mathcal{R}(F))$ contains a subdivision of $H$. Therefore, for every bijection $\pi:[|\mathcal{R}(F)|]\rightarrow \mathcal{R}(F)$, we have $\tilde{\mathcal{F}}_H\subseteq \mathcal{C}_{\sigma,F,\mathcal{R}(F)}\subseteq \mathcal{C}_{\sigma,F,\mathcal{R}(F),\pi}$. It follows that for every integer $t\geq 1$, we have $\tilde{\mathcal{F}}^t_H\subseteq \mathcal{C}_{\sigma,F,\mathcal{R}(F),\pi}$. This, together with Theorem~\ref{mainconnectify} and Lemma~\ref{cleanlemma}, implies Theorem~\ref{mainsubforestrestate}.
\end{proof}
In fact, one may deduce Theorem~\ref{mainstarforest1} directly using the material from Sections~\ref{daviessec}, \ref{noblocksec}, \ref{distancesec} and \ref{ramseyblocksec} (and in particular, skipping Section~\ref{connectedsec}).

\section{Obtaining a connectifier}\label{connectedsec} 
We begin with the following folklore result, see, for example, \cite{wallpaper} for a proof.

\begin{theorem}\label{minimalconnected}
Let $G$ be a connected graph, $X\subseteq V(G)$ with $|X|=3$ and $H$ be a connected induced subgraph of $G$ with $X\subseteq H$ and with $H$ minimal subject to inclusion. Then one of the following holds.
\begin{itemize}
\item There exists a vertex $a \in H$ and three paths $\{P_x:x\in X\}$ (possibly of length zero) where $P_x$ has ends $a$ and $x$, such that
\begin{itemize}
\item $H = \bigcup_{x\in X} P_x$, and;
\item the sets $\{P_x\setminus \{a\}: x\in X\}$ are pairwise disjoint and anticomplete.
\end{itemize}

\item There exists a triangle with vertex set $\{a_x:x\in X\}$ in $H$ and three paths $\{P_x: x\in X$\} (possibly of length zero) where $P_x$ has ends $a_x$ and $x$, such that
\begin{itemize}
\item $H= \bigcup_{x\in X} P_x$;
\item the sets $\{P_x \setminus \{a\}: x\in X\}$ are pairwise disjoint and anticomplete, and;
\item for distinct $x,y\in X$, $a_xa_y$ is the only edge of $H$ between $P_x$ and $P_y$.
\end{itemize}
\end{itemize}
\end{theorem}

Theorem~\ref{minimalconnected} may be reformulated as follows: for every choice of three vertices $x,y,z$ in a connected graph $G$, there is an induced subgraph $H$ of $G$ containing $x,y,z$ such that, for some $\delta\in [3]$, $H$ is isomorphic to either a subdivision of $K_{1,\delta}$ or the line graph of a subdivision of $K_{1,\delta}$, and $\mathcal{Z}(H)\subseteq \{x,y,z\}$. The main result of this section, the following, can be viewed as a qualitative extension of Theorem~\ref{minimalconnected}.

\begin{theorem}\label{minimalconnectedgeneral}
  For every integer $\eta\geq 1$, there exists an integer $\mu=\mu(\eta)\geq 1$ with the following property. Let $G$ be a graph and
  $S \subseteq V(G)$  with $|S|\geq \mu$ such that $S$ is contained in a connected component of $G$. Then $G$ contains an $(S,\eta)$-connectifier $H$.
  In particular,  $H$ is connected, $|H\cap S|=\eta$, and every vertex in $H\cap S$ has degree at most $\eta$ in $H$.
  \end{theorem}
  
For a graph $G$, $S\subseteq G$ and positive integer $\eta$, one may observe that $(S,\eta)$-connectifiers are minimal with respect to being connected and containing $\eta$ vertices from $S$. Also, for $\eta_1,\eta_2\geq 4$ (which, given Theorem~\ref{minimalconnected}, captures the main content of Theorem~\ref{minimalconnectedgeneral}) and distinct $i_1,i_2\in \{0,1,\ldots,4\}$, no $(S,\eta_1)$-connectifier of type $i_1$ contains an induced subgraph which is an $(S,\eta_2)$-connectifier of type $i_2$. Therefore, Theorem~\ref{minimalconnectedgeneral} provides an efficient characterization of all minimally connected induced subgraphs of $G$ containing many vertices from a sufficiently large subset $S$ of vertices in $G$.

In order to prove Theorem~\ref{minimalconnectedgeneral}, we need a few definitions and a result from \cite{Davies}. By a \textit{big clique} in a graph $J$, we mean a maximal clique of cardinality at least three. A graph $J$ is said to be a \textit{bloated tree} if
\begin{itemize}
    \item every edge of $J$ is contained in at most one big clique of $J$.
\item for every big clique $K$ of $J$ and every $v\in K$, $v$ has at most one neighbor in $J\setminus K$; and
\item the graph obtained from $J$ by contracting each big clique into a vertex is a tree.
\end{itemize}

It follows that every bloated tree is connected, and every connected induced subgraph of a bloated tree is a bloated tree. Furthermore, we deduce:

\begin{lemma}\label{blaotedlemma}
    Let $J$ be a bloated tree. Then for every cycle $C$ in $J$, $V(C)$ is a clique of $J$. 
\end{lemma}
\begin{proof}
    Suppose for a contradiction that for some cycle $C$ in $J$, $V(C)$ contains two vertices which are non-adjacent in $J$. Let $C$ be chosen with $|V(C)|=k$ as small as possible. It follows that $k\geq 4$. Let $C=c_1\dd \cdots \dd c_k\dd c_1$ such that $c_1$ and $c_i$ are not adjacent for some $i\in \{3,\ldots,k-1\}$. Let $P$ be a path in $J$ from $c_1$ to $c_i$ with $P^*\subseteq \{c_2\ldots, c_{i-1}\}$ and let $Q$ be a path in $J$ from $c_1$ to $c_i$ with $Q^*\subseteq \{c_{i+1}\ldots, c_{k}\}$. So $P$ and $Q$ are internally vertex-disjoint and $|P|,|Q|\geq 3$. Also, $H=J[P\cup Q]$ is a connected induced subgraph of $J$, and so $H$ is a bloated tree. If $P^*$ is anticomplete to $Q^*$, then $H$ is cycle. But then the graph obtained from $H$ by contracting each big clique into a vertex is $H$ itself, which is not tree, a contradiction with $H$ being a bloated tree. It follows that there exists $p\in P^*$ and $q\in Q^*$ such that $pq\in E(J)$. Consequently, $C_1=c_1\dd P\dd p\dd q\dd Q\dd c_1$ and $C_2=c_i\dd P\dd p\dd q\dd Q\dd c_i$ are two cycles in $J$ with $|V(C_1)|,|V(C_2)|<|V(C)|$. Thus, by the choice of $C$, for each $i\in \{1,2\}$, $K_i=J[V(C_i)]$ is a clique of $J$. For each $i\in \{1,2\}$, let $K_i'$ be a maximal clique of $J$ containing $K_i$. Then we have $c_1\in K_1'$ and $c_i\in K_2'$, which implies that $K_1'$ and $K_2'$ are distinct. But now the edge $pq\in E(J)$ is contained in two maximal cliques of $J$, namely $K_1'$ and $K_2$, which violates $J$ being a bloated tree. This proves Lemma~\ref{blaotedlemma}.
\end{proof}
The following was proved in \cite{Davies}:
\begin{theorem}[Davies \cite{Davies}]\label{daviesbloated} For every integer $k\geq 1$, there exists an integer $f=f(k)$ such that if $G$ is a connected graph and $S\subseteq V(G)$ with $|S|\geq f(k)$, then $G$ has an induced subgraph $J$ which is a bloated tree and $|J\cap S|\geq k$.
\end{theorem}

We also need the following well-known result; see, for example, \cite{wallpaper} for a proof.

\begin{lemma}\label{degreeorpath}
    For all positive integers $d,q$, there exists a positive integer $N(d,q)$ such that for every connected graph $G$ on at least $N(d,q)$ vertices, either $G$ contains a vertex of degree at least $d$, or there is a path in $G$ with $q$ vertices.
\end{lemma}

For a graph $G$ and a set $S\subseteq G$, by an $S$-\textit{bump} we mean a vertex $v\in G\setminus S$ of degree two in $G$, say $N_G(v)=\{v_1,v_2\}$, such that $v_1v_2\notin E(G)$. Also, by \textit{suppressing} the $S$-bump $v$ we mean removing $v$ from $G$ and adding the edge $v_1v_2$ (hence, $G$ is a subdivision of the resulting graph). We are now ready to prove Theorem~\ref{minimalconnectedgeneral}, which we restate:

\setcounter{theorem}{1}
\begin{theorem}\label{minimalconnectedgeneral}
  For every integer $\eta\geq 1$, there exists an integer $\mu=\mu(\eta)\geq 1$ with the following property. Let $G$ be a graph and
  $S \subseteq V(G)$  with $|S|\geq \mu$ such that $S$ is contained in a connected component of $G$. Then $G$ contains an $(S,\eta)$-connectifier $H$.
  In particular,  $H$ is connected, $|H\cap S|=\eta$, and every vertex in $H\cap S$ has degree at most $\eta$ in $H$.
  \end{theorem}

  \begin{proof}
     Let $f(\cdot)$ be as in Theorem~\ref{daviesbloated}, and  $N(\cdot, \cdot)$ be as in Lemma~\ref{degreeorpath}. We choose 
     $$\mu=\mu(\eta)=f(\max\{N(\eta,8\eta^2+\eta),2\}).$$
     
     By Theorem~\ref{daviesbloated}, since $|S|\geq \mu$, it follows that $G$ has an induced subgraph $J$ which is a bloated tree with $|J\cap S|\geq \max\{N(\eta,8\eta^2+\eta),2\}$, and subject to this property, $J$ has as few vertices as possible. Assume that $\eta=2$. Then, since $J$ is connected and $|J\cap S|\geq 2$, there is a path $H$ in $J$ with ends in $S$ and $H^*\cap S=\emptyset$. But then $H$ is an $(S,2)$-connectifier of type $2$ in $G$, as desired. Therefore, we may assume that $\eta\geq 3$.

     \sta{\label{getSminimal} Let $X\subseteq J$ such that $X$ is connected. Then for every connected component $Q$ of $J\setminus X$, we have $Q\cap S\neq \emptyset$. In particular, we have $\mathcal{Z}(J)\subseteq S$.}

     Suppose not. Let $Q$ be a component of $J\setminus X$ such that $Q\cap S= \emptyset$. Since $X$ connected, it holds that $J\setminus Q$ is connected, as well. It follows that $J\setminus Q$ is bloated tree and $|(J\setminus Q)\cap S|=|J\cap S|$, which contradicts the minimality of $J$. This proves \eqref{getSminimal}.\medskip 

     Let $J_1$ be the graph obtained from $J$ by successively suppressing $S$-bumps in $J$ until there are none. Then $J_1$ is also a bloated tree, and $J$ is a subdivision of $J_1$. The following is immediate from \eqref{getSminimal} and the definition of $J_1$.

      \sta{\label{getSminimalj1}$J_1$ has no $S$-bump and $J_1\cap S=J\cap S$. Also, for every $X\subseteq J_1$ with $X$ connected and every connected component $Q$ of $J_1\setminus X$, we have $Q\cap S\neq \emptyset$. In particular, we have $\mathcal{Z}(J_1)\subseteq S$.}

     Since $J$ is a bloated tree and so contains no hole, it follows that $J$ is a subdivision of $J_1$ with the additional property that for every edge $e\in E(J_1)$ which is contained in a big clique of $J_1$, we have $e\in E(J)$ (that is, $e$ is not subdivided while obtaining $J$ from $J_1$). This, along with the fact that $J_1\cap S=J\cap S$, implies that $J$ contains an $(S,\eta)$-connectifier if and only if $J_1$ contains an $(S,\eta)$-connectifier. Therefore, in order to prove Theorem~\ref{minimalconnectedgeneral}, it suffices to show that $J_1$ contains an $(S,\eta)$-connectifier, which we do in the rest of the proof.
           
      \sta{\label{maximaldisjoint} Let $K$ be a maximal clique of $J_1$, and for every $v\in K$, let $Q_v$ be the connected component of $J_1\setminus (K\setminus\{v\})$ containing $v$. Then for every two distinct vertices $u,v\in K$, we have $Q_u\cap Q_v=\emptyset$, and $uv$ is the only edge of $J_1$ between $Q_u$ and $Q_v$.}

    Suppose for a contradiction that there exist two distinct vertices $u,v\in K$ for which either $Q_u\cap Q_v\neq \emptyset$ or there is an edge in $J_1$ different from $uv$ with one end in $Q_u$ and one end in $Q_v$. It follows that $J_1[Q_u\cup Q_v]-uv$ is connected, and so there exists a path $P$ in $J_1$ of length more than one from $u$ to $v$ with $P^*\subseteq (Q_u\cup Q_v)\setminus \{u,v\}\subseteq J_1\setminus K$. Let $x\in P^*$. Then $C=u\dd P\dd v\dd u$ is a cycle in $J_1$. Since $J_1$ is a bloated tree, by Lemma~\ref{blaotedlemma}, $V(C)$ is a clique, and so $x$ is adjacent to both $u$ and $v$. Now, suppose that there exists a vertex $y\in K\setminus N_{J_1}(x)$. Then we have $y\notin \{u,v\}$, and so $C'=x\dd u\dd y\dd v\dd x$ is a cycle in $J_1$ where $V(C')$ contains two non-adjacent vertices, namely $x$ and $y$, which contradicts Lemma~\ref{blaotedlemma} and the fact that $J_1$ is a bloated tree. Therefore, $x$ is complete to $K$, and so $K\cup \{x\}$ is a clique of $J_1$ strictly containing $K$. This violates the maximality of $K$, and so proves \eqref{maximaldisjoint}.\medskip

    \sta{\label{largedegreehonest} Suppose that $J_1$ contains a big clique $K$ with $|K|\geq \eta$. Then $J_1$ contains an $(S,\eta)$-connectifier of type $0$.}

For every $v\in K$, let $Q_v$ be the connected component of $J_1\setminus (K\setminus\{v\})$ containing $v$. Then by \eqref{maximaldisjoint}, for every two distinct vertices $u,v\in K$, we have $Q_u\cap Q_v=\emptyset$, and there is no edge in $J_1$ with one end in $Q_u$ and one end in $Q_v$ except for $uv$. Also, by \eqref{getSminimalj1}, for every $v\in K$, we have $Q_v\cap S\neq \emptyset$. Therefore, since $Q_v$ is connected, we can choose a path $P_v$ in $Q_v$ from $v$ to a vertex $\ell_v\in S$ (possibly $v=\ell_v$) with $P_v\cap S=\{\ell_v\}$. It follows that for distinct $u,v\in K$, we have $P_u\cap P_v=\emptyset$, and there is no edge in $J_1$ with one end in $P_u$ and one end in $P_v$ except for $uv$. Now, let $K'\subseteq K$ with $|K'|=\eta$. Since $\eta\geq 3$, it follows that $H=J_1[\bigcup_{v\in K'}P_v]$ is a loosely $(S,\eta)$-tied line graph of a subdivided star; that is, $H$ is an $(S,\eta)$-connectifier of type $0$ in $J_1$. This proves \eqref{largedegreehonest}.\medskip

       \sta{\label{maximaldisjointvertex} Let $x\in J_1$ such that $N_{J_1}(x)$ is a stable set of $J_1$, and for every $a\in N_{J_1}(x)$, let $Q_a$ be the connected component of $J_1\setminus x$ containing $a$. Then the sets $\{Q_a: a\in N_{J_1}(x)\}$ are pairwise disjoint and anticomplete to each other.}

    Suppose for a contradiction that there exist two distinct vertices $a,b\in N_{J_1}(x)$ for which either $Q_a\cap Q_b\neq \emptyset$, or there is an edge in $J_1$ with one end in $Q_a$ and one end in $Q_b$. It follows that $J_1[Q_a\cup Q_b]$ is connected, and so there exists a path $P$ in $J_1$ of length more than one from $a$ to $b$ with $P^*\subseteq Q_a\cap Q_b\setminus \{a,b\}\subseteq J_1\setminus \{a,b,x\}$. Then $C=a\dd P\dd b\dd x\dd a$ is a cycle in $J_1$ where $V(C)$ contains two non-adjacent vertices, namely $a$ and $b$. This contradicts Lemma~\ref{blaotedlemma} and the fact that $J_1$ is a bloated tree, and so proves \eqref{maximaldisjointvertex}.\medskip
    
Now we can handle the case where $J_1$ contains vertices of large degree.

     \sta{\label{largedegree} Suppose that $J_1$ has a vertex of degree at least $\eta$. Then $J_1$ contains an $(S,\eta)$-connectifier of type $0$ or $1$.}

    Since $J_1$ is a bloated tree, for every vertex $x\in J_1$, either $N_{J_1}(x)$ is a clique, or $N_{J_1}(x)$ is stable set, or $J_1[N_{J_1}(x)]$ has an isolated vertex $y$ for which $N_{J_1}(x)\setminus \{y\}$ is a clique. Therefore, $J_1$ has a vertex of degree at least $\eta$, and it follows that either $J_1$ contains a big clique $K$ with $|K|\geq \eta$ or there exists a vertex $x\in V(J_1)$ of degree at least $\eta$ in $J_1$ such that $N_{J_1}(x)$ is a stable set of $J_1$. In the former case, \eqref{largedegree} follows from \eqref{largedegreehonest}. So we may assume that the latter case holds. For each $a\in N_{J_1}(x)$, let $Q_a$ be the connected component of $J_1\setminus x$ containing $a$. Then by \eqref{maximaldisjointvertex}, the sets $\{Q_a: a\in N_{J_1}(x)\}$ are pairwise disjoint and anticomplete to each other. Also, by \eqref{getSminimalj1}, for every $a\in N_{J_1}(x)$, we have $Q_a\cap S\neq \emptyset$. Therefore, since $Q_a$ is connected, we can choose a path $P_a$ in $Q_a$ from $a$ to a vertex $\ell_a\in S$ (possibly $a=\ell_a$) with $P_a\cap S=\{\ell_a\}$. It follows that the paths $\{P_a: a\in N_{J_1}(x)\}$ are pairwise disjoint and anticomplete to each other. Let $A$ be a subset of $N_{J_1}(x)$  with $|A|=\eta-1$ if $x\in S$ and $|A|=\eta$ if $x\notin S$. Then $H=J_1[\bigcup_{a\in A}P_a]$ is a $(S,\eta)$-tied rooted subdivided star with root $x$ such that $(H\cap S)\setminus \mathcal{L}(H)\subseteq \{x\}$; that is, $H$ is $(S,\eta)$-connectifier in $J_1$ of type $1$. This proves \eqref{largedegree}.\medskip

Henceforth, by \eqref{largedegree}, we may assume that $J_1$ has no vertex of degree at least $\eta$. Also, by \eqref{getSminimalj1},  we have $|J_1|\geq |J_1\cap S|\geq N(\eta,8\eta^2+\eta)$. As a result, by Lemma~\ref{degreeorpath}, $J_1$ contains a path $P$ on $8\eta^2+\eta$ vertices.

\sta{\label{pathcase} Suppose that there is no path in $P\setminus S$ of length $8\eta$. Then $J_1$ contains an $(S,\eta)$-connectifier of type $2$.}

Suppose not. Then $P$ contains no $(S,\eta)$-tied path. Let $|P\cap S|=s$. It follows that $s<\eta$. Therefore, since there is no path in $P\setminus S$ of length $8\eta$, we have $|P|\leq 8\eta(s+1)+s<8\eta^2+\eta$, a contradiction. This proves \eqref{pathcase}.\medskip 

In view of \eqref{pathcase}, we may assume that $P$ contains a path $P_1$ of length $8\eta$ with $P_1\cap S=\emptyset$, say
$$P_1=d_0\dd a_1\dd b_1\dd c_1\dd d_1\dd a_2\dd b_2\dd c_2\dd d_2\dd \cdots \dd a_{2\eta}\dd b_{2\eta}\dd c_{2\eta}\dd d_{2\eta}.$$
For each $i\in [2\eta]$, let $A_i=\{a_i,b_i,c_i\}$, let $L_i$ be the connected component of $J_1\setminus A_i$ containing $P_1[d_0,d_{i-1}]$, and let $R_i$ be the connected component of $J_1\setminus X_i$ containing $P_1[d_{i},d_{2\eta}]$. We deduce:

\sta{\label{disjointcomponentsxi} For each $i\in [2\eta]$, $L_i$ and $R_i$ are distinct, and so $L_i\cap R_i=\emptyset$.}

Suppose not. Then $J_1[L_i\cup R_i]$ is connected. Therefore, there exists a path $Z$ in $J_1$ from a vertex $z\in L_i$ to a vertex $z'\in R_i$ such that $Z^*\subseteq (L_i\cup R_i)\setminus P_1$. But then $C=z\dd P_1\dd z'\dd Z\dd z$ is a cycle in $J_1$ and $V(C)$ contains two non-adjacent vertices, namely $a_i$ and $c_i$, contradicting that $J_1$ is a bloated tree. This proves \eqref{disjointcomponentsxi}.

\sta{\label{thirdcomponent} For each $i\in [2\eta]$, there exists a component $Q_i$ of $J_1\setminus A_i$ different from $L_i$ and $R_i$.}

Suppose not. Then $J_1\setminus A_i$ has exactly two distinct components, namely $L_i$ and $R_i$. Assume that $b_i$ has degree two in $J_1$. Then, since $b_i\in P_1\subseteq J_1\setminus S$, it follows that $b_i$ is an $S$-bump, which violates \eqref{getSminimalj1}. So there exists a vertex $z\in N_{J_1}(b_i)\setminus A_i\subseteq L_i\cup R_i$, say $z\in L_i$. Consequently, since $L_i$ is connected, there exists a path $Z$ in $L_i$ from $z$ to a vertex $z'\in P_1[d_{0},d_{i-1}]$ with $Z\cap P_1=\{z'\}$. But then $C=b_i\dd z\dd Z\dd z'\dd P_1\dd b_i$ is a cycle in $J_1$ and $V(C)$ contains two non-adjacent vertices, namely $b_i$ and $d_{i-1}$, contradicting that $J_1$ is a bloated tree. This proves \eqref{thirdcomponent}.

\sta{\label{disjointcomponentsqi} For each $i\in [2\eta]$, let $Q_i$ be as in \eqref{thirdcomponent}. Then we have $P_1\cap Q_i=\emptyset$ and $N_{J_1}(Q_i)\subseteq A_i$. Also, the sets $\{Q_i: i\in [2\eta]\}$ are pairwise disjoint and anticomplete to each other.}

The first two assertions are immediate from the fact that  $Q_i$ is a component of $J_1\setminus A_i$ different from $L_i$ and $R_i$. For the third one, suppose for a contradiction that $Q_i\cup Q_j$ is connected for some distinct $i,j\in [2\eta]$, say $i<j$. Since $J_1$ is connected and $N_{J_1}(Q_j)\subseteq A_j$, it follows that $Q_j\cup A_j$ is connected, and so $Q_i\cup Q_j\cup A_j$ is connected. As a result, there exists a path $R$ in $J_1$ with one end $q\in Q_i$ and one end $q'\in A_j\subseteq R_i$ with $R^*\subseteq Q_j$. Also, we have $A_i\cap R\subseteq A_i\cap (Q_i\cup Q_j\cup A_j)\subseteq P_1\cap (Q_i\cap Q_j)=\emptyset$. In other words, $R$ is a path in  $J_1\setminus A_i$ from $q\in Q_i$ to $q'\in L_i$. But then we have $q\in Q_i\cap R_i$, a contradiction with \eqref{thirdcomponent}. This proves \eqref{disjointcomponentsqi}.\medskip

For each $i\in [2\eta]$, let $Q_i$ be as in \eqref{thirdcomponent}. Then by \eqref{getSminimalj1}, since $A_i$ is connected, we have $Q_i\cap S\neq \emptyset$. Also, from \eqref{thirdcomponent} and the connectivity of $J_1$, we have $N_{Q_i}(A_i)\neq \emptyset$. Therefore, since $Q_i$ is connected, we can choose a path $W_i$ in $Q_i$ from a vertex in $x_i\in N_{Q_i}(A_i)$ to a vertex in $y_i \in Q_i\cap S$ (possibly $x_i=y_i$) such that $W_i^*\cap (N_{Q_i}(A_i)\cup S)=\emptyset$.
Let $G_i=J_1[A_i\cup W_i]$. It follows that $G_i$ is connected and $G_i\cap S=\{y_i\}$.  

The following is easily observed:

\sta{\label{caterfinal}The sets $\{G_i:i\in [2\eta]\}$ are pairwise disjoint and anticomplete to each other. Also, for every $i\in [2\eta]$, $d_{i-1}a_i$ and $c_id_{i}$ are the only edges in $J_1$ with one end in $G_i$ and one end in $P_1\setminus G_i$.}

The proof is almost concluded. Note that since $J_1$ is a bloated tree, it follows that for every $i\in [2\eta]$, there is no cycle in $J_1$ containing both $a_i$ and $c_i$. Consequently, we have either $|N_{A_i}(x_i)|=1$, or $N_{A_i}(x_i)=\{a_i,b_i\}$, or $N_{A_i}(x_i)=\{b_i,c_i\}$, as otherwise $x_i\dd a_i\dd b_i\dd c_i\dd x_i$ is a cycle in $J_1$ containing both $a_i$ and $c_i$. Let $I\subseteq [2\eta]$. We say $I$ is \textit{light} if $|N_{A_i}(x_i)|=1$ for every $i\in I$. Also, we say $I$ is \textit{heavy} if for every $i\in I$, we have either $N_{A_i}(x_i)=\{a_i,b_i\}$, or $N_{A_i}(x_i)=\{b_i,c_i\}$. It follows that there exists $I\subseteq [2\eta]$ with $|I|=\eta$ which is either light or heavy. Let $i_1$ and $i_{\eta}$ be smallest and the largest elements of $I$, respectively.  It follows from $\eta\geq 3$ that $i_1$ and $i_{\eta}$ are distinct and $i_{\eta}\geq 3$. Let $Z_{1}$ be a path in $G_{i_1}$ from $c_{i_1}$ to $y_{i_1}$, and 
let $Z_{\eta}$ be a path in $G_{i_{\eta}}$ from $a_{i_{\eta}}$ to $y_{i_{\eta}}$. Let 
$$H=J_1\left[P_1[c_{i_1},a_{i_{\eta}}] \cup \left(Z_1\cup Z_{\eta}\right)\cup \left(\bigcup_{i\in I\setminus \{i_1,i_{\eta}\}}G_i\right)\right].$$
Using \eqref{caterfinal}, it is straightforward to observe that if $I$ is light, then $H$ is a loosely $(S,\eta)$-tied caterpillar, and if  $I$ is heavy, then $H$ is a loosely $(S,\eta)$-tied line graph of a caterpillar. In other words, $H$ is an $(S,\eta)$-connectifier of type $3$ or $4$. This completes the proof of Theorem~\ref{minimalconnectedgeneral}.
  \end{proof}
  
\setcounter{theorem}{0}

\section{Strong $k$-blocks}\label{noblocksec}
Let $G$ be a graph. By a \textit{separation} in $G$ we mean a triple $(L,M,R)$ of pairwise disjoint subsets of vertices in $G$ with $L\cup M\cup R=G$, such that neither $L$ nor $R$ is empty and $L$ is anticomplete to  $R$ in $G$. Let $x,y\in G$ be distinct. We say a set $M\subseteq G\setminus \{x,y\}$ \textit{separates $x$ and $y$} if there exists a separation $(L,M,R)$ in $G$ with $x\in L$ and $y\in R$. For a positive integer $k$, a $k$-\textit{block} in $G$ is a maximal set $B$ of at least $k$ vertices such that no two distinct vertices $x,y\in B$ are separated by a set $M\subseteq G\setminus \{x,y\}$ with $|M|<k$.  The application of $k$-blocks to bounding the treewidth in hereditary graph classes is not unprecedented; see for example, \cite{lozin,longhole2}. However, we find it best to work with a stronger notion of a $k$-block, which we define next.

Let $k$ be a positive integer and let $G$ be a graph. A \textit{strong $k$-block} in $G$ is a set $B$ of at least $k$ vertices in $G$ such that for every $2$-subset $\{x,y\}$ of $B$, there exists a collection $\mathcal{P}_{\{x,y\}}$ of at least $k$ distinct and pairwise internally disjoint paths in $G$ from $x$ to $y$, where for every two distinct $2$-subsets $\{x,y\}, \{x',y'\}\subseteq B$ and every choice of paths $P\in \mathcal{P}_{\{x,y\}}$ and $P'\in \mathcal{P}_{\{x',y'\}}$, we have $P\cap P'=\{x,y\}\cap \{x',y'\}$.

In this section, we prove that for all positive integers $k$ and $t$, every $t$-clean graph with no strong $k$-block has bounded treewidth. In other words, we show that for every positive integer $k$, the class of all graphs with no strong $k$-block is clean.

To begin with, we need some definitions as well as a couple of results from the literature. For a tree $T$ and an edge $xy\in E(T)$, we denote by $T_{x,y}$ the component of $T-{xy}$ containing $x$. Let $G$ be a graph and $(T,\chi)$ be a tree decomposition for $G$. For every $S\subseteq T$, let $\chi(S)=\bigcup_{x\in S}\chi(x)$. Also, for every edge $xy\in E(T)$, we define an \textit{adhesion} for $(T,\chi)$ as $\chi(x,y)=\chi(x)\cap \chi(y)=\chi(T_{x,y})\cap \chi(T_{y,x})$. For every $x\in V(T)$, by the \textit{torso at $x$}, denoted by $\hat{\chi}(x)$, we mean the graph obtained from the bag $\chi(x)$ by, for each $y\in N_T(x)$, adding an edge between every two non-adjacent vertices $u,v\in \chi(x,y)$. It is a well-known observation that clique cutsets do no effect the treewidth. More precisely, the following holds (a proof can be worked out easily using Lemma 5 from \cite{cliquetw}).

\begin{theorem}[folklore, see Lemma 5 in \cite{cliquetw}]\label{bodtorso}
Let $G$ be a graph and let $(T,\chi)$ be a tree decomposition for $G$. Then the treewidth of $G$ is at most the maximum treewidth of a torso $\hat{\chi}(x)$ taken over all $x\in V(T)$. 
\end{theorem}

Next we bring the material we need from \cite{tighttw} and \cite{Weissauerblock}. The \textit{fatness} of a tree decomposition $(T, \chi)$ of an $n$-vertex graph $G$ is the $(n + 1)$-tuple $(a_0,\ldots, a_n)$, where $a_i$ denotes the number of parts of $(T,\chi)$ of size $n-i$. If $(T, \chi)$
has lexicographically minimum fatness among all tree decompositions with all adhesions less
than $k$, we call $(T, \chi)$ \textit{$k$-atomic}. Also, a tree decomposition $(T,\chi)$ of a graph $G$ is \textit{tight} if for each vertex $x\in V(T)$ and every neighbor $y\in V(T)$ of $x$, there is a component $C$ of $\chi(T_{y,x})\setminus \chi(T_{x,y})$ such that every vertex in $\chi(x,y)$ has a neighbor in $C$. The following is proved in \cite{Weissauerblock}.

\begin{lemma}[Wei\ss auer, Lemma 6 in \cite{Weissauerblock}]\label{atomiclemma}
    Every $k$-atomic tree decomposition is tight.
\end{lemma}

Let $(T,\chi)$ be a tree decomposition for a graph $G$ and $\mathcal{S}$ be a set of pairwise disjoint subtrees of $T$. Let $T'$ be the tree obtained from $T$ by contracting every subtree $S\in \mathcal{S}$ into a new vertex $v_S$. Let $\chi':V(T')\rightarrow 2^{V(G)}$ be defined as follows. Let $\chi'(v_S)=\chi(S)$ for every $S\in \mathcal{S}$, and let $\chi'(v)=\chi(v)$ for every $v\in V(T')\setminus \{v_S:S\in \mathcal{S}\}=V(T)\setminus (\bigcup_{S\in \mathcal{S}}S)$. One may readily observe that $(T',\chi')$ is a tree decomposition for $G$, which is referred to as a \textit{contraction of} $(T,\chi)$. The following theorem from \cite{tighttw} is the key ingredient in our proof of the main result of this section.\footnote{We remark that the corresponding statement in \cite{tighttw}, namely ``Theorem 4'' therein, does not explicitly mention that $(T,\chi)$ is a contraction of a $k$-atomic tree decomposition. However, as the reader can check, the proof given in Section 3 of \cite{tighttw} is easily seen to yield this: it starts with a $k$-atomic tree decomposition ``$(T,\mathcal{V})$'' with $k=r(r-1)$, and concludes at the end that the desired tree decomposition is a certain contraction of $(T,\mathcal{V})$.}

\begin{theorem}[Erde and Wei\ss auer \cite{tighttw}, see also \cite{Grohe}]\label{tightdegorg}
Let $r$ be a positive integer, and let $G$ be a graph containing no subdivision of $K_r$ as a subgraph. Then $G$ admits a tree decomposition $(T,\chi)$ for which the following hold.
\begin{itemize}
\item $(T,\chi)$ is a contraction of a $k$-atomic tree decomposition for $G$ with $k=r(r-1)$.
\item Every adhesion of $(T,\chi)$ has cardinality less than $r^2$.
\item For every $x\in V(T)$, either $\hat{\chi}(x)$ has fewer than $r^2$ vertices of degree at least $2r^4$, or $\hat{\chi}(x)$ has no minor isomorphic to  $K_{2r^2}$.
\end{itemize}
\end{theorem}

It is straightforward to check that every contraction of a tight tree decomposition is tight. Also, for every positive integer $k$ and every graph $G$, if $G$ contains a subdivision of $K_{k^3}$ as a subgraph, then $G$ contains a strong $k$-block. Therefore, the following is immediate from Theorem~\ref{tightdegorg} and Lemma~\ref{atomiclemma}. 

\begin{figure}[t]
\centering

\begin{tikzpicture}[scale=2.5,auto=left]
\tikzstyle{every node}=[inner sep=1.2pt,circle,draw]  
\centering

\node (s10) at (0,1.2) {};
\node (s12) at (0.6,1.2){};
\node(s14)[fill]  at (1.2,1.2){};
\node(s16)[fill] at (1.8,1.2){};
\node(s18)[fill] at (2.4,1.2){};
\node(s110) at (3,1.2){};
\node(s112) at (3.6,1.2){};
\node(s114) at (4.2,1.2){};

\node (s20) at (0,0.9) {};
\node (s21) at (0.3,0.9) {};
\node(s22) at (0.6,0.9){};
\node (s23) at (0.9,0.9) {};
\node(s24)[fill] at (1.2,0.9){};
\node (s25)[fill] at (1.5,0.9) {};
\node(s26) [label=60:$1$] at (1.8,0.9){};
\node (s27)[fill] at (2.1,0.9) {};
\node(s28)[fill] at (2.4,0.9){};
\node (s29)[fill] at (2.7,0.9) {};
\node (s210)[fill] at (3,0.9) {};
\node (s211)[fill] at (3.3,0.9) {};
\node (s212) at (3.6,0.9) {};
\node (s213) at (3.9,0.9) {};
\node (s214) at (4.2,0.9) {};
\node (s215) at (4.5,0.9) {};

\node (s30) at (0,0.6) {};
\node (s31) at (0.3,0.6) {};
\node(s32) at (0.6,0.6){};
\node (s33) at (0.9,0.6) {};
\node(s34)[fill] at (1.2,0.6){};
\node (s35)[fill] at (1.5,0.6) {};
\node(s36)[fill] at (1.8,0.6){};
\node (s37)[fill] at (2.1,0.6) {};
\node(s38)[fill] at (2.4,0.6){};
\node (s39)[label=60:$2$] at (2.7,0.6) {};
\node (s310)[fill] at (3,0.6) {};
\node (s311)[fill] at (3.3,0.6) {};
\node (s312)[fill] at (3.6,0.6) {};
\node (s313)[fill] at (3.9,0.6) {};
\node (s314)[fill] at (4.2,0.6) {};
\node (s315) at (4.5,0.6) {};

\node (s40) at (0,0.3) {};
\node (s41)[fill] at (0.3,0.3) {};
\node(s42)[fill] at (0.6,0.3){};
\node (s43)[fill] at (0.9,0.3) {};
\node(s44)[fill] at (1.2,0.3){};
\node (s45)[fill] at (1.5,0.3) {};
\node(s46)[label=60:$4$] at (1.8,0.3){};
\node (s47)[fill] at (2.1,0.3) {};
\node(s48)[fill] at (2.4,0.3) {};
\node (s49)[fill] at (2.7,0.3) {};
\node (s410)[fill] at (3,0.3) {};
\node (s411)[fill] at (3.3,0.3) {};
\node (s412)[label=60:$3$] at (3.6,0.3) {};
\node (s413)[fill] at (3.9,0.3) {};
\node (s414)[fill] at (4.2,0.3) {};
\node (s415) at (4.5,0.3) {};

\node (s50) at (0,0) {};
\node (s51)[fill] at (0.3,0) {};
\node(s52)[fill] at (0.6,0){};
\node (s53)[label=60:$7$] at (0.9,0) {};
\node(s54)[fill] at (1.2,0){};
\node (s55)[fill] at (1.5,0) {};
\node(s56)[fill] at (1.8,0){};
\node (s57)[fill] at (2.1,0) {};
\node(s58)[fill] at (2.4,0) {};
\node (s59)[label=60:$5$] at (2.7,0) {};
\node (s510)[fill] at (3,0) {};
\node (s511)[fill] at (3.3,0) {};
\node (s512)[fill] at (3.6,0) {};
\node (s513)[fill] at (3.9,0) {};
\node (s514)[fill] at (4.2,0) {};
\node (s515) at (4.5,0) {};

\node (s60) at (0,-0.3) {};
\node (s61) at (0.3,-0.3) {};
\node(s62)[fill] at (0.6,-0.3){};
\node (s63)[fill] at (0.9,-0.3) {};
\node(s64)[fill] at (1.2,-0.3){};
\node (s65)[fill] at (1.5,-0.3) {};
\node(s66)[label=60:$8$] at (1.8,-0.3){};
\node (s67)[fill] at (2.1,-0.3) {};
\node(s68)[fill] at (2.4,-0.3) {};
\node (s69)[fill] at (2.7,-0.3) {};
\node (s610)[fill] at (3,-0.3) {};
\node (s611)[fill] at (3.3,-0.3) {};
\node (s612)[label=60:$6$] at (3.6,-0.3) {};
\node (s613)[fill] at (3.9,-0.3) {};
\node (s614)[fill] at (4.2,-0.3) {};
\node (s615) at (4.5,-0.3) {};

\node (s70) at (0,-0.6) {};
\node (s71) at (0.3,-0.6) {};
\node(s72) at (0.6,-0.6){};
\node (s73) at (0.9,-0.6) {};
\node(s74) at (1.2,-0.6){};
\node (s75)[fill] at (1.5,-0.6) {};
\node(s76)[fill] at (1.8,-0.6){};
\node (s77)[fill] at (2.1,-0.6) {};
\node(s78)[fill] at (2.4,-0.6) {};
\node (s79)[label=60:$9$] at (2.7,-0.6) {};
\node (s710)[fill] at (3,-0.6) {};
\node (s711)[fill] at (3.3,-0.6) {};
\node (s712)[fill] at (3.6,-0.6) {};
\node (s713)[fill] at (3.9,-0.6) {};
\node (s714) at (4.2,-0.6) {};
\node (s715) at (4.5,-0.6) {};

\node (s80) at (0,-0.9) {};
\node(s82) at (0.6,-0.9){};
\node(s84) at (1.2,-0.9){};
\node(s86) at (1.8,-0.9){};
\node(s88)[fill] at (2.4,-0.9) {};
\node (s810)[fill] at (3,-0.9) {};
\node (s812) at (3.6,-0.9) {};
\node (s814) at (4.2,-0.9) {};


\foreach \from/\to in {s10/s12, s12/s14,s18/s110,s110/s112,s112/s114}
\draw [dashed] (\from) -- (\to);

\foreach \from/\to in {s14/s16,s16/s18}
\draw [thick] (\from) -- (\to);

\foreach \from/\to in {s20/s21, s21/s22, s22/s23, s23/s24, s25/s26,s26/s27,s211/s212,s212/s213,s213/s214, s214/s215}
\draw [dashed] (\from) -- (\to);

\foreach \from/\to in {s24/s25,s27/s28,s28/s29,s29/s210,s210/s211}
\draw [thick] (\from) -- (\to);

\foreach \from/\to in {s30/s31, s31/s32, s32/s33, s33/s34, s38/s39, s39/s310, s314/s315}
\draw [dashed] (\from) -- (\to);

\foreach \from/\to in {s34/s35, s35/s36, s36/s37, s37/s38, s310/s311, s311/s312, s312/s313, s313/s314}
\draw [thick] (\from) -- (\to);

\foreach \from/\to in {s40/s41, s45/s46, s46/s47, s411/s412,s412/s413, s414/s415}
\draw [dashed] (\from) -- (\to);

\foreach \from/\to in {s41/s42, s42/s43, s43/s44, s44/s45, s47/s48, s48/s49, s49/s410, s410/s411, s413/s414}
\draw [thick] (\from) -- (\to);

\foreach \from/\to in {s50/s51, s52/s53, s53/s54, s58/s59,s59/s510, s514/s515}
\draw [dashed] (\from) -- (\to);

\foreach \from/\to in {s51/s52, s54/s55, s55/s56,s56/s57,s57/s58, s510/s511, s511/s512, s512/s513, s513/s514}
\draw [thick] (\from) -- (\to);

\foreach \from/\to in {s60/s61, s61/s62, s65/s66,s66/s67, s611/s612,s612/s613, s614/s615}
\draw [dashed] (\from) -- (\to);

\foreach \from/\to in {s62/s63, s63/s64, s64/s65, s67/s68, s68/s69, s69/s610, s610/s611, s613/s614}
\draw [thick] (\from) -- (\to);

\foreach \from/\to in {s70/s71, s71/s72, s72/s73, s73/s74, s74/s75, 
s78/s79,s79/s710, s713/s714, s714/s715}
\draw [dashed] (\from) -- (\to);

\foreach \from/\to in {s75/s76, s76/s77, s77/s78, s710/s711, s711/s712, s712/s713}
\draw [thick] (\from) -- (\to);

\foreach \from/\to in {s80/s82, s82/s84,s84/s86,s86/s88,s810/s812,s812/s814}
\draw [dashed] (\from) -- (\to);

\foreach \from/\to in {s88/s810}
\draw [thick] (\from) -- (\to);


\foreach \from/\to in {s10/s20, s30/s40, s50/s60, s70/s80}
\draw [dashed] (\from) -- (\to);

\foreach \from/\to in {s21/s31, s61/s71}
\draw [dashed] (\from) -- (\to);

\foreach \from/\to in {s41/s51}
\draw[thick] (\from) -- (\to);

\foreach \from/\to in {s12/s22, s32/s42, s72/s82}
\draw [dashed] (\from) -- (\to);

\foreach \from/\to in {s52/s62}
\draw[thick] (\from) -- (\to);

\foreach \from/\to in {s23/s33,s43/s53, s63/s73}
\draw [dashed] (\from) -- (\to);

\foreach \from/\to in {s74/s84}
\draw [dashed] (\from) -- (\to);

\foreach \from/\to in {s14/s24, s34/s44, s54/s64}
\draw[thick] (\from) -- (\to);

\foreach \from/\to in {s25/s35,s45/s55, s65/s75}
\draw [thick] (\from) -- (\to);

\foreach \from/\to in {s16/s26, s36/s46, s56/s66, s76/s86}
\draw [dashed] (\from) -- (\to);

\foreach \from/\to in {s27/s37,s47/s57, s67/s77}
\draw [thick] (\from) -- (\to);

\foreach \from/\to in {s18/s28, s38/s48, s58/s68, s78/s88}
\draw [thick] (\from) -- (\to);

\foreach \from/\to in {s29/s39,s49/s59, s69/s79}
\draw [dashed] (\from) -- (\to);

\foreach \from/\to in {s110/s210}
\draw [dashed] (\from) -- (\to);

\foreach \from/\to in {s310/s410, s510/s610, s710/s810}
\draw [thick] (\from) -- (\to);

\foreach \from/\to in {s211/s311,s411/s511, s611/s711}
\draw [thick] (\from) -- (\to);

\foreach \from/\to in {s112/s212, s312/s412, s512/s612, s712/s812}
\draw [dashed] (\from) -- (\to);

\foreach \from/\to in {s213/s313}
\draw [dashed] (\from) -- (\to);

\foreach \from/\to in {s413/s513, s613/s713}
\draw [thick] (\from) -- (\to);

\foreach \from/\to in {s114/s214, s714/s814}
\draw [dashed] (\from) -- (\to);

\foreach \from/\to in {s314/s414, s514/s614}
\draw [thick] (\from) -- (\to);

\foreach \from/\to in {s215/s315,s415/s515, s615/s715}
\draw [dashed] (\from) -- (\to);

\draw[thick] (1.8,1.2) circle (0.05);

\draw[thick] (1.2,0.9) circle (0.05);
\draw[thick] (2.4,0.9) circle (0.05);

\draw[thick] (1.5,0.6)  circle (0.05);
\draw[thick] (2.1,0.6)  circle (0.05);
\draw[thick] (3.3,0.6)  circle (0.05);

\draw[thick] (0.6,0.3) circle (0.05);
\draw[thick] (1.2,0.3) circle (0.05);
\draw[thick] (2.4,0.3) circle (0.05);
\draw[thick] (3,0.3) circle (0.05);
\draw[thick] (4.2,0.3) circle (0.05);

\draw[thick] (0.3,0) circle (0.05);
\draw[thick] (1.5,0) circle (0.05);
\draw[thick] (2.1,0) circle (0.05);
\draw[thick] (3.3,0) circle (0.05);
\draw[thick] (3.9,0) circle (0.05);

\draw[thick] (1.2,-0.3) circle (0.05);
\draw[thick] (2.4,-0.3) circle (0.05);
\draw[thick] (3,-0.3) circle (0.05);
\draw[thick] (4.2,-0.3) circle (0.05);

\draw[thick] (2.1,-0.6) circle (0.05);
\draw[thick] (3.3,-0.6) circle (0.05);
\draw[thick] (3.9,-0.6) circle (0.05);

\draw[thick] (3,-0.9) circle (0.05);


\tikzstyle{every node}=[inner sep=1.2pt, fill=black,circle]  
\centering

\node(t10)[draw] at (1.2,-1.4) {};
\draw[thick] (1.2,-1.4) circle (0.05);

\node[draw] at (1.5,-1.4) {};

\node[inner sep=0, outer sep=0, fill=white] at (1.5,-1.55) {$1$};


\node(t12)[draw] at (1.8,-1.4){};
\draw[thick] (1.8,-1.4) circle (0.05);

\node[draw] at (1.95,-1.4){};
\node[draw] at (2.1,-1.4){};
\node[draw] at (2.25,-1.4){};

\node[inner sep=0, outer sep=0, fill=white] at (2.1,-1.55) {$2$};

\node(t14)[draw] at (2.4,-1.4){};
\draw[thick] (2.4,-1.4) circle (0.05);

\node[draw] at (2.55,-1.4){};
\node[draw] at (2.7,-1.4){};
\node[draw] at (2.85,-1.4){};

\node[inner sep=0, outer sep=0, fill=white] at (2.7,-1.55) {$3$};

\node(t16)[draw] at (3,-1.4){};
\draw[thick] (3,-1.4) circle (0.05);


\node(t20)[draw] at (1.2,-1.7) {};
\draw[thick] (1.2,-1.7) circle (0.05);

\node[draw] at (1.2,-1.55) {};
\node[draw] at (1.35,-1.7) {};


\node(t21)[draw] at (1.5,-1.7) {};
\draw[thick] (1.5,-1.7) circle (0.05);

\node[draw] at (1.65,-1.7) {};
\node[draw] at (1.5,-1.85) {};


\node(t22)[draw] at (1.8,-1.7){};
\draw[thick] (1.8,-1.7) circle (0.05);

\node[draw] at (1.95,-1.7) {};
\node[draw] at (1.8,-1.55) {};

\node[inner sep=0, outer sep=0, fill=white] at (1.8,-1.85) {$4$};


\node(t23)[draw] at (2.1,-1.7){};
\draw[thick] (2.1,-1.7) circle (0.05);

\node[draw] at (2.25,-1.7) {};
\node[draw] at (2.1,-1.85) {};


\node(t24)[draw] at (2.4,-1.7){};
\draw[thick] (2.4,-1.7) circle (0.05);

\node[draw] at (2.4,-1.55) {};
\node[draw] at (2.55,-1.7) {};

\node[inner sep=0, outer sep=0, fill=white] at (2.4,-1.85) {$5$};


\node(t25)[draw] at (2.7,-1.7){};
\draw[thick] (2.7,-1.7) circle (0.05);

\node[draw] at (2.7,-1.85) {};
\node[draw] at (2.85,-1.7) {};


\node(t26)[draw] at (3,-1.7){};
\draw[thick]  (3,-1.7) circle (0.05);

\node[draw] at  (3,-1.55) {};
\node[draw] at  (3.15,-1.7) {};

\node[inner sep=0, outer sep=0, fill=white] at (3,-1.85) {$6$};


\node(t27)[draw] at (3.3,-1.7){};
\draw[thick] (3.3,-1.7) circle (0.05);

\node[draw] at (3.3,-1.85) {};


\node(t30)[draw] at (1.2,-2) {};
\draw[thick] (1.2,-2) circle (0.05);

\node[draw] at (1.35,-2) {};

\node(t31)[draw] at (1.5,-2) {};
\draw[thick] (1.5,-2) circle (0.05);

\node[draw] at (1.65,-2) {};

\node[inner sep=0, outer sep=0, fill=white] at (1.5,-2.15) {$7$};

\node(t32)[draw] at (1.8,-2){};
\draw[thick]  (1.8,-2) circle (0.05);

\node[draw] at  (1.8,-2.15) {};
\node[draw] at  (1.95,-2) {};

\node(t33)[draw] at (2.1,-2){};
\draw[thick] (2.1,-2) circle (0.05);

\node[draw] at (2.25,-2) {};

\node[inner sep=0, outer sep=0, fill=white] at (2.1,-2.15) {$8$};

\node(t34)[draw] at (2.4,-2){};
\draw[thick]  (2.4,-2) circle (0.05);

\node[draw] at  (2.4,-2.15) {};
\node[draw] at  (2.55,-2) {};

\node(t35)[draw] at (2.7,-2){};
\draw[thick] (2.7,-2) circle (0.05);

\node[draw] at (2.85,-2) {};

\node[inner sep=0, outer sep=0, fill=white] at (2.7,-2.15) {$9$};

\node(t36)[draw] at (3,-2){};
\draw[thick]  (3,-2) circle (0.05);

\node[draw] at  (3,-2.15) {};
\node[draw] at  (3.15,-2) {};


\node(t37)[draw] at (3.3,-2){};
\draw[thick]  (3.3,-2) circle (0.05);


\node(t40)[draw] at (1.2,-2.3) {};
\draw[thick] (1.2,-2.3) circle (0.05);

\node[draw] at  (1.2,-2.15) {};
\node[draw] at  (1.35,-2.3) {};
\node[draw] at  (1.5,-2.3) {};
\node[draw] at  (1.65,-2.3) {};

\node(t42)[draw] at (1.8,-2.3){};
\draw[thick] (1.8,-2.3) circle (0.05);

\node[draw] at  (1.95,-2.3) {};
\node[draw] at  (2.1,-2.3) {};
\node[draw] at  (2.25,-2.3) {};

\node(t44)[draw] at (2.4,-2.3){};
\draw[thick] (2.4,-2.3) circle (0.05);

\node[draw] at  (2.6,-2.3) {};
\node[draw] at  (2.8,-2.3) {};


\node(t46)[draw] at (3,-2.3){};
\draw[thick] (3,-2.3) circle (0.05);


\foreach \from/\to in {t10/t12, t12/t14,t14/t16}
\draw [thick] (\from) -- (\to);

\foreach \from/\to in {t20/t21, t21/t22, t22/t23, t23/t24, t24/t25, t25/t26,t26/t27}
\draw [thick] (\from) -- (\to);

\foreach \from/\to in {t30/t31, t31/t32, t32/t33, t33/t34, t34/t35, t35/t36, t36/t37}
\draw [thick] (\from) -- (\to);

\foreach \from/\to in {t40/t42, t42/t44, t44/t46}
\draw [thick] (\from) -- (\to);

\foreach \from/\to in {t10/t20, t30/t40}
\draw [thick] (\from) -- (\to);

\foreach \from/\to in {t21/t31}
\draw [thick] (\from) -- (\to);

\foreach \from/\to in {t12/t22, t32/t42}
\draw [thick] (\from) -- (\to);

\foreach \from/\to in {t23/t33}
\draw [thick] (\from) -- (\to);

\foreach \from/\to in {t14/t24, t34/t44}
\draw [thick] (\from) -- (\to);

\foreach \from/\to in {t25/t35}
\draw [thick] (\from) -- (\to);

\foreach \from/\to in {t16/t26,t36/t46}
\draw [thick] (\from) -- (\to);

\foreach \from/\to in {t27/t37}
\draw [thick] (\from) -- (\to);

\end{tikzpicture}

\caption{Proof of \eqref{wallintorso}: the subgraph of $W_{8 \times 8}$ induced by the filled nodes (top) is isomorphic to a proper subdivision of $W_{4\times 4}$ (bottom). Note the correspondence between the numbered removed nodes at the top and the inner faces of the subdivided wall at the bottom.}
\label{fig:gridin5x5wall}
\end{figure}
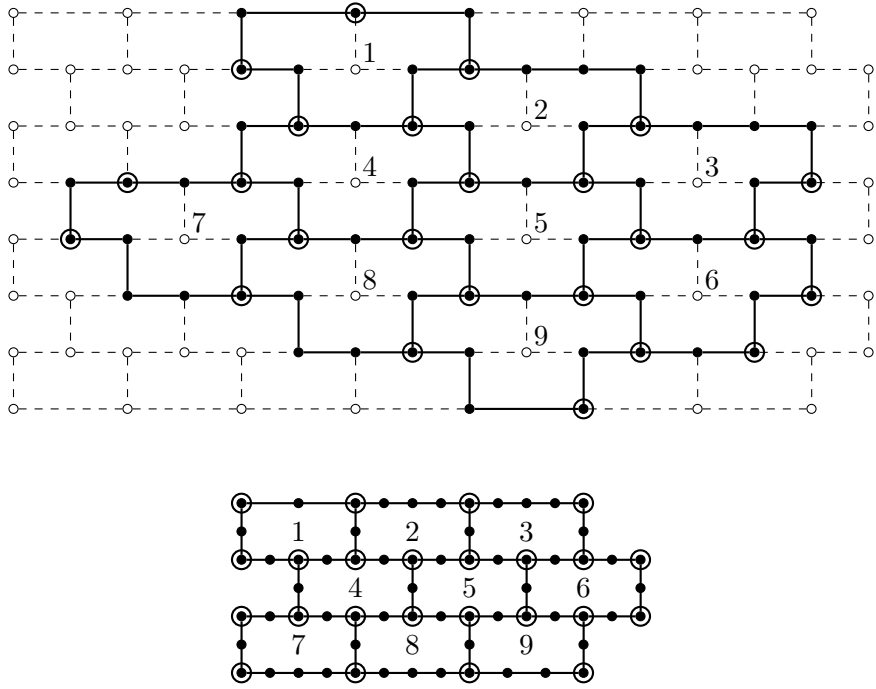

\begin{theorem}\label{tightdeg}
Let $k$ be a positive integer and let $G$ be a graph containing no strong $k$-block. Then $G$ admits a tight tree decomposition $(T,\chi)$ for which the following hold.
\begin{itemize}
\item Every adhesion of $(T,\chi)$ has cardinality less than $k^6$.
\item For every $x\in V(T)$, either $\hat{\chi}(x)$ has fewer than $k^6$ vertices of degree at least $2k^{12}$, or $\hat{\chi}(x)$ has no minor isomorphic to $K_{2k^6}$.
\end{itemize}
\end{theorem}
We can now prove the main result of this section. For every positive integer $k$, let $\mathcal{B}_k$ be the class of all graphs with no strong $k$-block.
\begin{theorem}\label{noblock}
For every integer $k\geq 1$, the class $\mathcal{B}_k$ is clean.
\end{theorem}
\begin{proof}

Let $t\geq 1$ and let $G\in \mathcal{B}_{k}^t$, that is, $G$ is a $t$-clean graph with no strong $k$-block. We aim to show that there exists an integer $w(k,t)\geq 1$ such that $\tw(G)\leq w(k,t)$. By Theorem~\ref{tightdeg}, $G$ has a tight tree decomposition $(T,\chi)$ for which every torso either has fewer than $k^6$ vertices of degree at least
$2k^{12}$ or has no minor isomorphic to $K_{2k^6}$. For each $x\in V(T)$, let $K_x\subseteq \hat{\chi}(x)$ be the set of all vertices in $\hat{\chi}(x)$ of degree at least $2k^{12}$. We define $\tau_x$ as follows: if $|K_x|<k^6$, then let $\tau_x=\hat{\chi}(x)\setminus K_x$, and otherwise let $\tau_x=\hat{\chi}(x)$. It follows that either $\tau_x$ has maximum degree less than $2k^{12}$, or $\tau_x$ has no minor isomorphic to $K_{2k^6}$. Let $\xi(\cdot,\cdot)$ be as in Theorem~\ref{minorclosed} and $\gamma(\cdot,\cdot)$ be as in Theorem~\ref{boundeddeg}. Let 
$$\gamma_0=\gamma(3,t),$$ 
$$\gamma_1=\gamma(2k^{12},2\gamma_0),$$
$$\xi_1=\xi(K_{2k^6},2\gamma_0),$$
$$w_1=w_1(k,t)=\max\{\gamma_1,\xi_1\}.$$
We claim that:

\sta{\label{wallintorso} For every $x\in V(T)$, we have $\tw(\tau_x)\leq w_1$.}

Suppose for a contradiction that $\tw(\tau_x)>w_1$ for some $x\in V(T)$. Note that either $\tau_x$ has maximum degree less than $2k^{12}$ or $\tau_x$ has no minor isomorphic to $K_{2k^6}$. Therefore, the choice of $w_1$ together with Theorems~\ref{minorclosed} and \ref{boundeddeg} implies that $\tau_x$ contains an induced subgraph $W$ which is isomorphic to either a subdivision of $W_{2\gamma_0\times 2\gamma_0}$ or the line graph of a subdivision of $W_{2\gamma_0\times 2\gamma_0}$. On the other hand, it can be seen that for every positive integer $q$, every subdivision of $W_{2q\times 2q}$ contains an induced subgraph isomorphic to a proper subdivision of $W_{q\times q}$ (see Figure~\ref{fig:gridin5x5wall}.) Consequently, $W$, and so $\tau_x$, contains an induced subgraph $W_0$ which is isomorphic to either a proper subdivision of $W_{\gamma_0\times \gamma_0}$ or the line graph of a proper subdivision of $W_{\gamma_0\times \gamma_0}$. In particular, $W_0$ has maximum degree at most three. Let us say a non-empty subset $K\subseteq W_0$ is a \textit{blossom} if there exists $y\in N_T(x)$ such that $K\subseteq \chi(x,y)$, and subject to this property, $K$ is maximal with respect to inclusion. It follows that every blossom $K$ is a clique in $W_0$ and so we have $|K|\in \{1,2,3\}$. Also, every two blossoms intersect in at most one vertex, and since no two triangles in $W_0$ share a vertex, blossoms of cardinality three are pairwise disjoint. Let $\mathcal{K}$ be the set of all blossoms, and for every blossom $K\in \mathcal{K}$, let us fix $y_K\in N_T(x)$ such that $K\subseteq \chi(x,y_K)$. From the maximality of blossoms, it follows that the vertices $\{y_K: K\in \mathcal{K}\}$ are all distinct. Note that $(T,\chi)$ is tight, and so for every $y\in N_T(x)$, there exists a component $C(y)$ of $\chi(T_{y,x})\setminus \chi(T_{x,y})$ such that the every vertex in $\chi(x,y)$ has a neighbor in $C(y)$. Since $(T,\chi)$ is a tree decomposition, it follows that the sets $\{C(y_K): y\in N_T(x)\}$ are pairwise distinct, disjoint and anticomplete in $G$. Let $H_K$ be a connected induced subgraph of $G[(C(y_K)\cup K)]$ which contains $K$, and subject to this property, assume that $H_K$ is minimal with respect to inclusion. It follows that if $|K|=1$, then $H_K=K$, if $|K|=2$, then $H_K$ is a path in $G$ between the two vertices in $K$ with $H_K^*\subseteq C(y_{K})$, and if $|K|=3$, then $H_K$ satisfies one of the two outcomes of Theorem~\ref{minimalconnected}. Also, the sets $\{H_K\setminus K: K\in \mathcal{K}\}$ are pairwise distinct, disjoint and anticomplete in $G$. Now, let 
\[H=G\left[\left(W_0\setminus \left(\bigcup_{K\in \mathcal{K}}K\right)\right)\cup \left(\bigcup_{K\in \mathcal{K}}H_K\right)\right].\]
Let $H'$ be the minor of $H$ obtained through the following steps in order:
\begin{enumerate}[(i)]
    \item For every blossom $K\in \mathcal{K}$ with $|K|=3$, contract the connected induced subgraph $H_K$ of $H$ into a vertex.
    \item For every blossom $K\in \mathcal{K}$ with $|K|=2$ such that $K$ is contained in a triangle of $W_0$, contract the path $H_K$ in $H$ into an edge between the two vertices in $K$.
    \item Contract each triangle of the resulting graph after (ii) into a vertex.
\end{enumerate}
Since $W_0$ is isomorphic to either a proper subdivision of $W_{\gamma_0\times \gamma_0}$ or the line graph of a proper subdivision of $W_{\gamma_0\times \gamma_0}$, it is readily observed that $H'$ is isomorphic to a subdivision of $W_{\gamma_0\times \gamma_0}$. It follows that $H$ contains $W_{\gamma_0\times \gamma_0}$ as a minor, and so we have $\tw(H)\geq \gamma_0=\gamma(3,t)+1$. Note that since $W_0$ has maximum degree at most three, $H$ has maximum degree at most three, as well. Therefore, by Theorem~\ref{boundeddeg}, $H$, and so $G$, contains either a subdivision of $W_{t\times t}$ or the line graph of a subdivision of $W_{t\times t}$ as a induced subgraph. But this violates the assumption that $G$ is $t$-clean, and so proves \eqref{wallintorso}.\medskip

    Now,  for every $x\in V(T)$, if $|K_x|<k^6$, then we have $\hat{\chi}(x)=\tau_x\cup K_x$, and otherwise we have $\hat{\chi}(x)=\tau_x$. This, along with \eqref{wallintorso}, implies that $\tw(\hat{\chi}(x))\leq w_1+k^6$ for every $x\in V(T)$. Hence, writing $w(k,t)=w_1(k,t)+k^6$, by Theorem~\ref{bodtorso}, we have $\tw(G)\leq w(k,t)$. This completes the proof of Theorem~\ref{noblock}.
    \end{proof}

Note that for every integer $k\geq 1$, if a graph $G$ contains a strong $k$-block, then $G$ contains $K_k$ as a topological minor, which in turn implies that $G$ contains every $k$-vertex graph as a topological minor. Therefore, the following common strengthening of Theorems~\ref{minorclosed} and \ref{boundeddeg} is in fact an immediate corollary of Theorem~\ref{noblock}:
\begin{corollary}
    For every graph $H$, the class of all graphs with no $H$-topological-minor is clean.
\end{corollary}

\section{$k$-blocks with distant vertices}\label{distancesec}
The main result of this section, Theorem~\ref{distance}, asserts that for every positive integer $k$, every graph containing a sufficiently large block contains either a subgraph that is a subdivision of a large complete graph with all paths short, or an induced subgraph which contains a $k$-block with its vertices pairwise far from each other. This will be of essential use in subsequent sections, and before proving it, we recall the classical result of Ramsey (see e.g.\ \cite{ajtai} for an explicit bound).
\begin{theorem}[See \cite{ajtai}]\label{classicalramsey}
For all integers $a,b\geq 1$, there exists an integer $R=R(a, b)\geq 1$ such that every graph $G$ on at least $R(a,b)$ vertices contains either a clique of cardinality $a$ or a stable set of cardinality $b$. In particular, for all integers $t\geq 1$ and $\rho\geq R(t,t)$, every graph $G$ containing $K_{\rho,\rho}$ as a subgraph contains either $K_t$ or $K_{t,t}$ as an induced subgraph.
\end{theorem}

For a graph $G$ and a positive integer $d$, a \textit{$d$-stable set} in $G$ is a set $S\subseteq G$ such that for every two distinct vertices $u,v\in S$, there is no path of length at most $d$ in $G$ from $u$ to $v$. Note that a $d$-stable set is also a $d'$-stable set for every $0<d'\leq d$. Here comes the main result of this section.

\begin{theorem}\label{distance}
    For all integers $d,k\geq 1$ and $m\geq 2$, there exists an integer $k_0=k_0(d,k,m)\geq 1$ with the following property. Let $G$ be a graph and $B_0$ be a strong $k_0$-block in $G$. Assume that $G$ does not contain a $(\leq d)$-subdivision of $K_m$ as a subgraph. Then there exists $A\subseteq G$ with $S\subseteq B_0\setminus A$ such that $S$ is both a strong $k$-block and a $d$-stable set in $G\setminus A$.
\end{theorem}
\begin{proof}
    Let $R(m,k)$ be as in Theorem~\ref{classicalramsey}. We show that $$k_0=k_0(d,k,m)=\binom{R(m,k)}{2}(d-1)+R(m,k)$$
    satisfies Theorem~\ref{distance}. Let $X\subseteq B_0$ with $|X|=R(m,k)$. Let $g=\binom{R(m,k)}{2}$. Let $e_1, \ldots, e_{g}$ be an enumeration of all $2$-subsets of $X$, and let $e_i=\{x_i,y_i\}$ for each $i\in [g]$.  Let $U_0=\emptyset$, and for every $i\in [g]$, having defined $U_{i-1}$, we define $P_i$ and $U_i$ as follows. If there exists a path $P$ in $G$ of length at most $d$ from $x_i$ to $y_i$ with $P^*\cap (U_{i-1}\cup X)=\emptyset$, then let $P_i=P$ and $U_i=U_{i-1}\cup P^*_i$. Otherwise, let $P_i=\emptyset$ and $U_i=U_{i-1}$. It follows that for all $i,j\in [g]$ with $i<j$ and $P_i,P_j\neq \emptyset$, we have $P_i\cap P^*_j=U_i\cap P^*_j=\emptyset$ and $P^*_i\cap P_j=P_i^*\cap X=\emptyset$.
    
    Let $G_0$ be the graph with $V(G_0)=X$ and for each $i\in [g]$, $x_i$ is adjacent to $y_i$ in $G_0$ if and only if $P_i\neq \emptyset$.
    
    \sta{\label{distancenoclique}$G_0$ contains no clique of cardinality $m$.}

    Suppose for a contradiction that $G_0$ contains a clique $C$ of cardinality $m$. Then for every $i\in [g]$ with $e_i\subseteq C$, we have $P_i\neq \emptyset$. Also, for all distinct $i,j\in [g]$, we have $P_i\cap P^*_j=P^*_i\cap P_j=\emptyset$. But then $G[\bigcup_{e_i\subseteq C}P_i]$, and so $G$ contains a ($\leq d$)-subdivision of $K_m$ as a subgraph, a contradiction. This proves \eqref{distancenoclique}.\medskip
    
 Since $|G_0|=|X|=R(m,k)$, it follows from Theorem~\ref{classicalramsey} and \eqref{distancenoclique} that $G_0$ contains a stable set $S$ of cardinality $k$. Let $A=U_g\cup (X\setminus S)$. Then we have $|A|\leq g(d-1)+R(m,k)-k$. Therefore, since $S\subseteq B_0\setminus A$ and $B_0$ is a strong $(g(d-1)+R(m,k))$-block, we deduce that $S$ is a strong $k$-block in $G\setminus A$.  It remains to show that $S$ is a $d$-stable set in $G\setminus A$. Suppose not. Then there exists $x,y\in S$ and a path $Q$ in $G\setminus A$ of length at most $d$ from $x$ to $y$. Thus, we may choose $i\in [g]$ such that $e_i\subseteq Q\cap S$. Therefore, assuming $P=Q[x_i,y_i]$, we have $P^*\cap S=\emptyset$. Now $P$ is a path in $G\setminus A$ (and so in $G$) of length at most $d$ from $x_i$ to $y_i$ with $P^*\subseteq G\setminus (A\cup S)=G\setminus (U_g\cup X)\subseteq G\setminus (U_{i-1}\cup X)$. It follows that $P_i\neq\emptyset$. But we have $e_i\subseteq S$ and $S$ is a stable set in $G_0$, which implies that $P_i=\emptyset$, a contradiction.
This completes the proof of Theorem~\ref{distance}.
\end{proof}

\section{Planted subdivided star forests}\label{ramseyblocksec}
In this section we extend ideas from \cite{lozin} to produce a subdivided star forest
whose roots are contained in sets with useful properties. Let $G$ be a graph, $S\subseteq G$, and $F$ a subdivided star forest. We say a subgraph $F'$ of $G$ isomorphic to $F$ is $S$-\textit{planted} if $F'$ is rooted and $\mathcal{R}(F')\subseteq S$. Write $\mathcal{H}_\lambda$ for the class of graphs with no holes of length greater than $\lambda$. The main result of this section is the following.

\begin{theorem}\label{blockmanystar}
For all positive integers $d,k,t,\delta,\lambda$, and  $\theta$ with $\delta\geq 2$, there exists a positive integer $k_1=k_1(d,k,t,\delta,\lambda, \theta)$ with the following property. Let $G$ be a $t$-clean graph and let $B_1$ be a strong $k_1$-block in $G$. Then there exist $A\subseteq V(G)$ and $S\subseteq B_1\setminus A$ such that the following hold.
\begin{itemize}
\item $S$ is both a strong $k$-block and a $d$-stable set in $G\setminus A$.
\item $G\setminus A$ contains an $S$-planted copy of $\theta S_{\delta,\lambda}$.
\item $G\setminus A$ contains a hole of length greater than $\lambda$.
\end{itemize}
In particular, we have $\mathcal{F}^t_{\theta S_{\delta,\lambda}},\mathcal{H}^t_{\lambda}\subseteq \mathcal{B}_{k_1}$.
\end{theorem}

Note that Theorem~\ref{blockmanystar}, combined with Theorem~\ref{noblock} and Lemma~\ref{cleanlemma}, implies Theorems~\ref{longholepili} and \ref{mainstarforest1} at once. Theorem~\ref{blockmanystar} is also a key tool in the proof of Theorem~\ref{mainconnectify} in Section~\ref{connectifysec}. We need the following two results from \cite{lozin}.

\begin{lemma}[Lozin and Razgon \cite{lozin}]\label{ramsey1}
For all positive integers $a$ and $b$, there is a positive integer $c=c(a, b)$ such that if a graph
$G$ contains a collection of $c$ pairwise disjoint subsets of $V(G)$, each of cardinality at most $a$ and with at least
one edge between every two of them, then $G$ contains $K_{b,b}$ as a subgraph.
\end{lemma}

\begin{theorem}[Lozin and Razgon \cite{lozin}]\label{ramsey2}
For all positive integers $p$ and $r$, there exists a positive integer $m = m(p,r)$ such that every
graph $G$ containing a ($\leq p$)-subdivision of $K_m$ as a subgraph contains either $K_{p,p}$ as a subgraph or a
proper ($\leq  p$)-subdivision of $K_{r,r}$ as an induced subgraph.
\end{theorem} 

We deduce the following lemma.

\begin{lemma} \label{lem:noKmsubdivision}
    For every integer $t\geq 1$, there exists an integer $n=n(t)\geq 1$ with the following property. Let $G$ be a $t$-clean graph and let $\rho$ be an integer with $\rho\geq R(t,t)$, where $R(\cdot,\cdot)$ is as in Theorem~\ref{classicalramsey}. Then $G$ does not contain a $(\leq \rho)$-subdivision of $K_n$ as a subgraph.
\end{lemma}

\begin{proof}
    Let $n=n(t)= m(R(t, t), 2t^2)$, where $m(\cdot, \cdot)$ is as in Theorem~\ref{ramsey2}. Suppose for a contradiction that $G$ contains a $(\leq \rho)$-subdivision of $K_n$ as a subgraph. Then by Theorem~\ref{ramsey2}, $G$ either contains $K_{\rho, \rho}$ as a subgraph, or contains an induced subgraph $H$ isomorphic to a proper subdivision of $K_{2t^2,2t^2}$. In the former case, by Theorem~\ref{classicalramsey}, $G$ contains either $K_t$ or $K_{t,t}$, which violates the assumption that $G$ is $t$-clean. In the latter case, note that a proper subdivision of $K_{2t^2,2t^2}$ contains a proper subdivision of every bipartite graph on at most $2t^2$ vertices. In particular, $H$, and so $G$, contains a subdivision of $W_{t\times t}$, again contradicting that $G$ is $t$-clean. This proves the Lemma~\ref{lem:noKmsubdivision}.
\end{proof}

We are now ready to prove the main result of this section.

\begin{proof}[Proof of Theorem~\ref{blockmanystar}]
    Let $R(\cdot,\cdot)$ be as in Theorem~\ref{classicalramsey}. Let $c = c(\lambda, R(t, t))$, where $c(\cdot, \cdot)$ is as in Lemma~\ref{ramsey1}. Let  $n = n(t)$, be as in Lemma~\ref{lem:noKmsubdivision}.  Let $k_0(\cdot,\cdot,\cdot)$ be as in Theorem~\ref{distance}. Let
    $$k_1 = k_1(d,k,t,\delta,\lambda, \theta)= k_0(\max\{d, R(t,t), 2\lambda + 1\}, \max\{k, R(c, \delta), \theta\}, n).$$
    We claim that this choice of $k_1$ satisfies Theorem~\ref{blockmanystar}. To see this, suppose that $G$ is a $t$-clean graph which has a strong $k_1$-block $B_1$. Note first that, by Lemma~\ref{lem:noKmsubdivision}, $G$ does not contain a $(\leq \max\{d, R(t,t), 2\lambda + 1\})$-subdivision of $K_n$ as a subgraph. Therefore, by Theorem~\ref{distance}, there exist $A \subseteq G$ and $S \subseteq B_1 \setminus A$ such that $S$ is both a strong $\max\{k, R(c, \delta), \theta\}$-block and a $\max\{d, R(t,t), 2\lambda + 1\}$-stable set in $G \setminus A$. In particular, $S$ is both a strong $k$-block and a $d$-stable set in $G \setminus A$, which proves the first bullet of Theorem~\ref{blockmanystar}. Next we claim that:

    \sta{\label{1star}For every $x\in S$, there exists a copy $F_x$ of $S_{\delta,\lambda}$ in $G\setminus A$ where $x\in F_x$ has degree $\delta$ in $F_x$.}
    
    It is easily seen that $|S|\geq 2$. Pick a vertex $y \in S\setminus \{x\}$. Since $S$ is a strong $R(c, \delta)$-block in $G\setminus A$, there exists a collection $\{P_i: i\in [R(c, \delta)]\}$ of pairwise internally disjoint paths in $G \setminus A$ from $x$ to $y$. Since $S$ is a $(2\lambda + 1)$-stable set in $G \setminus A$, for each $i\in [R(c, \delta)]$, $P_i$ has length greater than $\lambda + 1$. Let $P_i'$ be the subpath of $P_i$ of length $\lambda$ containing $x$ as an end. Then $\{P_i':i\in [R(c,\delta)]\}$ is a collection of $R(c, \delta)$ pairwise disjoint subsets of $G\setminus A$, each of cardinality $\lambda$. 
    Let $\Gamma$ be the graph with $V(\Gamma)=[R(c, \delta)]$ such that for all distinct $i,j\in [R(c, \delta)]$, $i$ is adjacent to $j$ in $\Gamma$ if and only if $P_i'\setminus \{x\}$ is not anticomplete to $P_j'\setminus \{x\}$ in $G$. By Theorem~\ref{classicalramsey}, $\Gamma$ contains either a clique of cardinality $c$ or a stable set of cardinality $\delta$. Suppose first that $\Gamma$ contains a clique of cardinality $c$. Then Lemma~\ref{ramsey1} implies that $G$ contains $K_{R(t, t), R(t, t)}$ as a subgraph, and thus by Theorem~\ref{classicalramsey}, $G$ contains $K_t$ or $K_{t, t}$, which violates the assumption that $G$ is $t$-clean. Consequently, $\Gamma$ has a stable set $I$ of cardinality $\delta$. But now $F_x=G[\bigcup_{i\in I}P_i']$ is a copy of $S_{\delta,\lambda}$ in $G\setminus A$ where $x\in F_x$ has degree $\delta$ in $F_x$. This proves \eqref{1star}.\medskip
    
Now we can prove the second bullet of Theorem~\ref{blockmanystar}. For every $x\in S$, let $F_x$ be as in \eqref{1star}. Note that since $S$ is a $(2\lambda + 1)$-stable set in $G \setminus A$, it follows that for all distinct $x,x'\in S$, $F_x$ and $F_{x'}$ are disjoint and anticomplete to each other. Also, since $S$ is a strong $\theta$-block, there exists $S'\subseteq S$ with $|S'|=\theta$. But now $G[\bigcup_{x\in S'}F_x]$ is an $S$-planted copy of $\theta S_{\delta,\lambda}$ in $G\setminus A$, as desired.

    It remains to prove the third bullet of Theorem~\ref{blockmanystar}. Proceeding as in the proof of \eqref{1star}, we choose distinct vertices $x,y\in S$ and two internally disjoint paths $P_1$ and $P_2$ in $G\setminus A$ from $x$ to $y$ such that $P_1'\setminus \{x\}$ is anticomplete to $P_2'\setminus \{x\}$, where for each $i\in \{1,2\}$, $P_i'$ is the subpath of $P_i$ of length $\lambda$ containing $x$ as an end. Traversing $P_1$ from $x$ to $y$, let $z$ be the first vertex in $P_1^*$ with a neighbor in $P_2 \setminus \{x\}$ (this vertex exists, since the neighbor of $y$ in $P_1$ is adjacent to $P_2 \setminus \{x\}$). Also, traversing $P_2$ from $x$ to $y$, let $w \in P_2\setminus \{x\}$ be the first neighbor of $z$ in $P_2\setminus \{x\}$. Note that since $P_1'$ is anticomplete to $P_2'$, it follows that either $z \notin P'_1$ or $w\notin P'_2$.  But now $x \dd P_1 \dd z \dd w \dd P_2 \dd x$ is a hole in $G\setminus A$ of length at least $\lambda+3$. This completes the proof of Theorem~\ref{blockmanystar}.
\end{proof}

\section{Proof of Theorem~\ref{mainconnectify}}\label{connectifysec}

The last step in the proof of Theorem~\ref{mainconnectify} is the following. Note that the condition $\delta\geq 3$ is due to the fact that there is only one choice of roots for subdivided star forests in which every component has a branch vertex, and so it is slightly more convenient to work with them. 

\begin{lemma}\label{blockconnectify}
    For all positive integers $t, \delta, \lambda, \sigma, \theta$ with $\delta\geq 3$ and $\theta\geq 2$, there exists an integer $k_2=k_2(t, \delta, \lambda, \sigma, \theta)\geq 1$ with the following property. Let $G$ be a $t$-clean graph containing a strong $k_2$-block. Then $G$ contains a $\sigma$-connectification of $(\theta S_{\delta,\lambda},\mathcal{R}(\theta S_{\delta,\lambda}))$. In other words, we have $\mathcal{C}^t_{\sigma,\theta S_{\delta,\lambda}, \mathcal{R}(\theta S_{\delta,\lambda})}\subseteq \mathcal{B}_{k_2}$.
\end{lemma}
\begin{proof}
    Let $\mu(\cdot)$ be as in Theorem~\ref{minimalconnectedgeneral}. Let
    $$\gamma_1=\mu(\max\{t,\sigma\theta,\theta+1\}),$$
    $$\gamma_2=\mu(\gamma_1),$$
    $$\gamma_3=\gamma_2((2t\gamma_1+\delta)\lambda+1).$$
    Let $k_1(\cdot, \cdot,\cdot,\cdot,\cdot, \cdot)$ be as in Theorem~\ref{blockmanystar}. We define:
    $$\displaystyle k_2=k_2(t, \delta, \lambda, \sigma, \theta)=k_1\left(2\sigma-1,\gamma_3+R(t,t)\binom{\gamma_3}{2t},t,2t\gamma_1+\delta,\lambda,\gamma_2\right).$$
    Let $B_2$ be a strong $k_2$-block in $G$. By Theorem~\ref{blockmanystar}, there exist $A\subseteq G$ and $S\subseteq B_2\setminus A$ such that the following hold. Let $G_0=G\setminus A$.
\begin{itemize}
\item $S$ is both a strong  $R(t,t)\binom{\gamma_3}{2t}$-block and a $(2\sigma -1)$-stable set in $G_0$.
    \item $G_0$ contains an $S$-planted copy $F$ of $\gamma_2 S_{2t\gamma_1+\delta,\lambda}$.
\end{itemize}

      Then $|\mathcal{R}(F)|=\gamma_2$ and $|F|=\gamma_3$. For every $x\in \mathcal{R}(F)$, let $F_x$ be the component of $F$ with root $x$. Let $W$ be the set of all vertices in $G_0\setminus F$ with at least $2t$ neighbors in $F$.

    \sta{\label{boundednbrs} We have $\displaystyle |W|<R(t,t)\binom{\gamma_3}{2t}$.}

    Suppose not. Let $q=R(t,t)\binom{\gamma_3}{2t}$ and let $w_1, \ldots, w_q\in W$ be distinct. For every $i\in [q]$, let $N_i$ be a set of $2t$ neighbors of $w_i$ in $F$. It follows that there exist $I\subseteq [q]$ and $N\subseteq F$ such that $|I|=R(t,t)$, $|N|=2t$ and $N_{i}=N$ for all $i\in I$. Note that since $F$ is a forest, $N$ contains a stable set $N'$ of $G_0$ with $|N'|=t$. Also, since $G_0$ is $t$-clean, it does not contains a clique of cardinality $t$. Thus, by Lemma~\ref{classicalramsey}, $G_0[\{w_i:i\in I\}]$ contains a stable set $N''$ of cardinality $t$. But then $G_0[N'\cup N'']$ is isomorphic to $K_{t,t}$, which contradicts that $G_0$ is $t$-clean. This proves \eqref{boundednbrs}.\medskip

    Let $G_1=G_0\setminus W$. Then $G_1$ is a $t$-clean induced subgraph of $G$. In order to prove Theorem~\ref{blockconnectify}, it suffices to show that $G_1$ contains a $\sigma$-connectification of $\theta S_{\delta,\lambda}$, which we do in the rest of the proof.
    
    Recall that $S$ is both a strong $(\gamma_3+R(t,t)\binom{\gamma_3}{2t})$-block and a $(2\sigma-1)$-stable set in $G_0$. Thus, since $S\setminus W\subseteq G_1$, by \eqref{boundednbrs}, $S\setminus W$ is both a strong $\gamma_3$-block and a $(2\sigma-1)$-stable set in $G_1$. Also, we have $\mathcal{R(F)}\subseteq S\setminus W$. It follows that $\mathcal{R}(F)$ is a $(2\sigma-1)$-stable set in $G_1$, and for every two distinct vertices $x,x'\in \mathcal{R}(F)$, since $|F\setminus \mathcal{R}(F)|<\gamma_3$, there is a path in $G_1\setminus (F\setminus \mathcal{R}(F))$ from $x$ to $x'$. Consequently, $G_1\setminus (F\cup \mathcal{R}(F))$ has a component containing $\mathcal{R}(F)$. Let $G_2$ be the graph obtained from $G_1$ by contracting $F_x$ into $x$ for each $x\in \mathcal{R}(F)$. Then $G_2$ contains $G_1\setminus (F\cup \mathcal{R}(F))$ as a spanning subgraph, and so $G_2$ has a component containing $\mathcal{R}(F)$. Since $\mathcal{R}(F)\geq \gamma_2=\mu(\gamma_1)$, from Theorem ~\ref{minimalconnectedgeneral} applied to $G_2$ and $\mathcal{R}(F)$, it follows that $G_2$ contains a connected induced subgraph $H_2$ such that, assuming $S'=H_2\cap \mathcal{R}(F)$, we have $|S'|= \gamma_1$ and every vertex in $S'$ has degree at most $\gamma_1$ in $H_2$. Let
$$H_1=G_1\left[H_2\cup \left(\bigcup_{x\in S'} F_x\right)\right].$$
In other words, $H_1$ is the induced subgraph of $G_1$ obtained from $H_2$ by undoing the contraction of $F_x$ into $x$ for each $x\in H_2\cap \mathcal{R}(F)$. It follows that $H_1$ is a connected induced subgraph of $G_1$ and $H_1\cap \mathcal{R}(F)=H_2\cap \mathcal{R}(F)=S'$. Moreover, since $\mathcal{R(F)}$ is a $(2\sigma-1)$-stable set in $G_1$, $S'$ is also a $(2\sigma-1)$-stable set in $H_1$.

\sta{\label{degree2t} For every $x\in S'$, we have $|N_{F_x}(H_1\setminus F_x)|< 2t\gamma_1$.}

Note that $N_{H_1\setminus F_x}(F_x)=N_{H_2}(x)$, and so $|N_{H_1\setminus F_x}(F_x)|\leq \gamma_1$. Also, since $H_1$ is an induced subgraph of $G_1$, by the definition of $W$, no vertex in $N_{H_1\setminus F_x}(F_x)\subseteq G_1\setminus F$ has at least $2t$ neighbors in $F_x$. Therefore, we have $|N_{F_x}(H_1\setminus F_x)|< 2t\gamma_1$. This proves \eqref{degree2t}.\medskip

The following is immediate from \eqref{degree2t} and the fact that for every $x\in S'$, $F_x$ is isomorphic to $S_{2t\gamma_1+\delta,\lambda}$.

\sta{\label{anticompletestars} For every $x\in S'$, $F_x$ contains an induced copy $F'_x$ of $S_{\delta,\lambda}$ containing $x$ such that $F'_x\setminus \{x\}$ is anticomplete to $H_1\setminus F'_x$.}

Next, we define:
$$H'_1=H_1\setminus \left(\bigcup_{x\in S'} (F'_x\setminus \{x\})\right).$$
It follows that $H_1'$ is a connected induced subgraph of $G_1$ and $S'\subseteq H'_1$ is a $(2\sigma-1)$-stable set in $H_1'$.

\sta{\label{getconnectifier}
$H_1'$, and so $G_1$, contains an $(S',\theta, \sigma)$-connectifier $H$ of type $i$ for some $i\in [4]$.}

Since $|S'|\geq \gamma_1=\mu(\max\{t,\theta \sigma,\theta+1\})$, we can apply Theorem~\ref{minimalconnectedgeneral} to $H_1'$ and $S'$. It follows that $H_1'$ contains an $(S',\max\{t,\theta \sigma,\theta+1\})$-connectifier $H'$. Since $S'$ is a $(2\sigma-1)$-stable set in $H_1'$, $H'\cap S'$ is also a $(2\sigma-1)$-stable set in $H'$. It is straightforward to observe that if $H'$ is of type $i$ for $i\in \{2,3,4\}$, then $H'$, and so $H_1'$, contains an $(S',\theta, \sigma)$-connectifier $H$. Also, if $H'$ is of type $0$, then $H'$ contains a clique of cardinality $t$, which violates that $G_1$ is $t$-clean. It remains to consider the case where $H'$ is of type $1$. Then $H'$ contains an $(S',\theta+1)$-tied rooted subdivided star $H''$ with root $r$ in which every stem has length at least $\sigma$ and $(H''\cap S')\setminus \mathcal{L}(H'')\subseteq \{r\}$.  Since $\theta\geq 2$, it follows that $H''$ has at least three vertices and $r$ is not a leaf of $H''$. If $H''$ is a path with ends $h_1,h_2\in S'$, then $\theta=2$ and $r\in S'$. This, along with the fact that $H''\cap S'$ is a $(2\sigma-1)$-stable set in $H''$, implies that $H=H''[h_1,r]$ has length at least $2\sigma$. But then $H$ is a $(S',\theta,\sigma)$-connectifier of type $2$ in $H''$, and so in $H_1'$. Also, if $H''$ is not a path, then $r$ is the unique branch vertex of $H''$. Again, since $H''\cap S'$ is $(2\sigma-1)$-stable set in $H'$ (and so in $H''$), there exists a stem $P$ of $H''$ such that every stem of $H''$ other than $P$ has length at least $\sigma$. Therefore, $H=H'\setminus (P\setminus \{r\})$ is an  $(S',\theta,\sigma)$-connectifier of type $1$ in $H''$, and so in $H_1'$. This proves \eqref{getconnectifier}.\medskip

Let $H$ be as in \eqref{getconnectifier}. Let $X=H\cap S$. Let $F'=\bigcup_{x\in X}F'_x$ and $\Xi=G_1[H\cup F']$. Then by \eqref{anticompletestars}, $F'$ is an induced subgraph of $\Xi$ isomorphic to $\theta S_{\delta,\lambda}$ and $F'\setminus X$ is anticomplete to $\Xi\setminus F$. Also, we have $\Xi\setminus (F'\setminus X)=H$. But then by \eqref{getconnectifier}, $\Xi$ is a $\sigma$-connectification of $(F',X)$, and so $\Xi$ is an induced subgraph of $G$ isomorphic to a $\sigma$-connectification of $(\theta S_{\delta,\lambda},\mathcal{R}(\theta S_{\delta,\lambda}))$. This completes the proof of Lemma~\ref{blockconnectify}.
\end{proof}

We need one more definition before proving Theorem~\ref{mainconnectify}. For two rooted subdivided star forests $F_1$ and $F_2$, we say $F_2$ \textit{embeds in} $F_1$ if $\mathcal{R}(F_2)\subseteq \mathcal{R}(F_1)$ and there exists a collection $\mathcal{S}$ of stems of $F_1$ such that $F_2=F_1\setminus ((\bigcup_{P\in \mathcal{S}}P)\setminus \mathcal{R}(F_1))$. 

Now we prove Theorem~\ref{mainconnectify}, which we restate:

\setcounter{section}{4}
\setcounter{theorem}{0}

\begin{theorem}
    Let $\sigma\geq 1$ be an integer, let $F$ be a rooted subdivided star forest of size at least two and let $\pi:[|\mathcal{R}(F)|]\rightarrow \mathcal{R}(F)$ be a bijection. Then the class $\mathcal{C}_{\sigma, F, \mathcal{R}(F),\pi}$ is clean. 
\end{theorem}
\begin{proof}
    Let $F$ be of maximum degree $\delta\geq 0$, reach $\lambda\geq 0$ and size $\theta\geq 2$. For every $x\in \mathcal{R}(F)$, let $F_x$ be the component of $F$ with root $x$. Let $F^+=\theta S_{\delta+3,\lambda+1}$ be rooted (with its unique choice of roots). For every $y\in \mathcal{R}(F^+)$, let $F^+_y$ be the component of $F^+$ with root $y$. Then for every $x\in \mathcal{R}(F)$ and every $y\in \mathcal{R}(F^+)$, $F^+_{y}$ contains a copy $F^+_{x,y}$ of $F_{x}$ such that $F^+_{x,y}$ embeds in $F^+_{y}$. Now, for every choice of bijections $\pi:[\theta]\rightarrow \mathcal{R}(F)$ and $\pi^+:[\theta]\rightarrow \mathcal{R}(F^+)$, and every $\sigma$-connectification $\Xi^+$ of $(F^+,\mathcal{R}(F^+))$ with respect to $\pi^+$, let $$\Xi=(\Xi^+\setminus F^+)\cup \left(\bigcup_{i\in [\theta]}F^+_{\pi(i),\pi^+(i)}\right).$$
    It follows that $\Xi$ is isomorphic to a $\sigma$-connectification of $(F,\mathcal{R}(F))$ with respect to $\pi$. In other words, for every bijection $\pi:[\theta]\rightarrow \mathcal{R}(F)$, every $\sigma$-connectification of $(F^+,\mathcal{R}(F^+))$ contains an induced subgraph isomorphic to a $\sigma$-connectification of $(F,\mathcal{R}(F))$ with respect to $\pi$. Therefore, we have $\mathcal{C}_{\sigma, F, \mathcal{R}(F),\pi}\subseteq \mathcal{C}_{\sigma, F^+, \mathcal{R}(F^+)}$. This, together with Lemma~\ref{blockconnectify}, implies that for every integer $t\geq 1$, we have $\mathcal{C}^t_{\sigma, F, \mathcal{R}(F),\pi}\subseteq \mathcal{C}^t_{\sigma, F^+, \mathcal{R}(F^+)}\subseteq  \mathcal{B}_{k_2}$, where $k_2=k_2(t,\delta+3,\lambda+1,\sigma,\theta)$ is as in Lemma~\ref{blockconnectify}. Now the result follows from Theorem~\ref{noblock} and Lemma~\ref{cleanlemma}.
\end{proof}

\setcounter{section}{9}
\setcounter{theorem}{0}

\section{Acknowledgement}
We thank Rose McCarty for bringing Theorem~\ref{tightdegorg} to our attention, and Daniel Wei\ss auer for helpful comments.

\bibliographystyle{plain}

\begin{thebibliography}{30}

  \bibitem{aboulker}
 P. Aboulker, I. Adler, E. J. Kim, N. L. D. Sintiari, and N. Trotignon.
\newblock {``On the treewidth of even-hole-free graphs.''}
{\em European Journal of Combinatorics} {\bf 98}, (2021), 103394. 


\bibitem{wallpaper}
  T. Abrishami, M. Chudnovsky, C. Dibek, S. Hajebi, P. Rz\k{a}\.{z}ewski, S. Spirkl, and K. Vu\v{s}kovi\'c. ``Induced subgraphs and tree decompositions II. Toward walls and their line graphs in graphs of bounded degree.''
  {\em arXiv:2108.01162}, (2021). 


\bibitem{pyramiddiamond}
T. Abrishami, M. Chudnovsky, S. Hajebi, and S. Spirkl. ``Induced subgraphs and tree decompositions IV. 
(Even hole, diamond, pyramid)-free graphs.'' {\em arXiv:2203.06775}, (2022).  

\bibitem{onenbr}
T. Abrishami, B. Alecu, M. Chudnovsky, S. Hajebi, and S. Spirkl. ``Induced subgraphs and tree decompositions V. One neighbor in a hole'' {\em arXiv:2205.04420}, (2022).  

\bibitem{ajtai}
 M. Ajtai, J. Koml\'{o}s and E. Szemer\'{e}di. \newblock {``A note on Ramsey numbers.''}
{\em  J. Combinatorial Theory, Ser. A} {\bf 29}, (1980), 354–360.

\bibitem{tw13}
B. Alecu, M. Chudnovsky, S. Hajebi, and S. Spirkl. ``Induced subgraphs and tree decompositions XIII. Basic obstructions in $\mathcal{H}$-free graphs for finite $\mathcal{H}$,'' {\em manuscript} (2023).  

\bibitem{multistar}
 M. Bonamy, \'{E}. Bonnet, H. D\'{e}pr\'{e}s, L. Esperet, C. Geniet, C. Hilaire, S. Thomass\'{e} and A. Wesolek, ``Sparse graphs with bounded induced cycle packing number have logarithmic treewidth'',
 {\em arXiv:2206.00594}, (2022).

 \bibitem{cliquetw}
H. Bodlaender and A. Koster. ``Safe separators for treewidth.'' {\em Discrete Mathematics} {\bf 306}, 3 (2006), 337--350.

\bibitem{Bodlaendersurvey}
H. L. Bodlaender. Treewidth: Structure and Algorithms. {\em SIROCCO 2007: Proceedings of the 14th International Colloquium on Structural Information and Communication Complexity}, (2007), 11--25.




\bibitem{daviesconst}
 J. Davies, appeared in an Oberwolfach technical report {\em DOI:10.4171/OWR/2022/1}.


\bibitem{Davies}
 J. Davies, ``Vertex-minor-closed classes are $\chi$-bounded.'' 
\emph{arXiv:2008.05069}, (2020).

\bibitem{tighttw}
 J. Erde and D. Wei\ss auer. \newblock {``A short derivation of the structure theorem for graphs with excluded topological minors.''}
{\em  SIAM Journal of Discrete Mathematics} {\bf 33}, 3 (2019), 1654--1661.

\bibitem{longholes}
P. Gartland, D. Lokshtanov, M. Pilipczuk, M. Pilipczuk, P. Rz\k{a}\.{z}ewski. ``Finding large induced sparse subgraphs in $C_{>t}$-free graphs in quasipolynomial time,'' {\em STOC 2021: Proceedings of the 53rd Annual ACM SIGACT Symposium on Theory of Computing}, (2021), 330--341.

\bibitem{Grohe}
M. Grohe and D. Marx. ``Structure theorem and isomorphism test for graphs with excluded topological subgraphs,'' {\em SIAM Journal on Computing} \textbf{44}, 1 (2015), 114--159.
		

\bibitem{Korhonen}
  T. Korhonen, ``Grid induced minor theorem for graphs of small degree.''
 {\em J. Combin. Theory Ser. B} \textbf{160} (2023), 206 -- 214.

\bibitem{lozin}
 V. Lozin, I. Razgon. \newblock {``Tree-width dichotomy.''}
{\em  European Journal of Combinatorics} {\bf 103}, (2022), 103517.

\bibitem{Menger}  K. Menger, ``Zur allgemeinen Kurventheorie.'' \textit {Fund. Math.}
      {\bf 10}, 1927, 96--115.
  

\bibitem{RS-GMII}
N. Robertson and P. Seymour. ``Graph minors. II. Algorithmic aspects of tree-width.'' \textit{J. of Algorithms}, {\bf 7} (3) (1986), 309–322.
   

\bibitem{RS-GMV}
N. Robertson and P. Seymour. ``Graph minors. V. Excluding a planar graph.'' \textit{J. Combin. Theory Ser. B}, {\bf 41} (1) (1986), 92–114.
   


  \bibitem{layered-wheels}
N.L.D. Sintiari and N. Trotignon. ``(Theta, triangle)-free and (even-hole, K4)-free graphs. Part 1: Layered
wheels,'' {\em J. Graph Theory} {\bf 97} (4) (2021), 475--509.

  \bibitem{longhole2}
 D. Wei\ss auer. \newblock {``In absence of long chordless cycles, large tree-width becomes a local phenomenon.''}
{\em  J. Combin. Theory Ser. B} {\bf 139}, (2019) 342–352.

  \bibitem{Weissauerblock}
 D. Wei\ss auer. \newblock {``On the block number of graphs.''}
{\em  SIAM Journal of Discrete Mathematics} {\bf 33}, (2019), 346–357.

\end{thebibliography}

\end{document}